%% file: manuscript.tex
\documentclass[openacc]{rsproca_new}

\def\Real{{\mathbb{R}}}

\title{Formation Mechanism of a Basin of Attraction for Passive Dynamic Walking Induced by Intrinsic Hyperbolicity}

\author{Ippei Obayashi$^{1}$, Shinya Aoi$^{2}$, Kazuo Tsuchiya$^{2}$, and
  Hiroshi Kokubu$^{3}$}

\address{
  $^{1}$ Advanced Institute for Materials Research (AIMR),
  Tohoku University,
  2--1--1 Katahira, Aoba-ku, Sendai 980--8577, Japan \\
  $^{2}$ Department of Aeronautics and Astronautics,
  Graduate School of Engineering,
  Kyoto University,
  Kyoto daigaku-Katsura, Nishikyo-ku, Kyoto,
  615--8540,
  Japan \\
  $^{3}$ Department of Mathematics, Graduate School of Science,
  Kyoto University,
  Kitashirakawa Oiwakecho, Sakyo-ku, Kyoto 606--8502,
  Japan 
}

\subject{Mechanical engeneering, Applied mathematics}

\corres{Ippei Obayashi\\
  \email{ippei.obayashi.d8@tohoku.ac.jp}}

\begin{document}

\begin{abstract}
Passive dynamic walking is a useful model for investigating the mechanical functions of the body that produce energy-efficient walking.
The basin of attraction is very small and thin, and it has a fractal-like shape; this explains the difficulty in producing stable passive dynamic walking.
The underlying mechanism that produces these geometric characteristics was not known.
In this paper, we consider this from the viewpoint of dynamical systems theory, and we use the simplest walking model to clarify the mechanism that forms the basin of attraction for passive dynamic walking.
We show that the intrinsic saddle-type hyperbolicity of the upright equilibrium point in the governing dynamics plays an important role in the geometrical characteristics of the basin of attraction; this contributes to our understanding of the stability mechanism of bipedal walking.
\end{abstract}

\maketitle

\section{Introduction}

When humans walk, the stance leg is almost straight, and it rotates around the contact point of the foot like an inverted pendulum.
Therefore, the center of mass (COM) is at its highest position during the midstance phase and at its lowest position during the double-support phase.
In contrast, the locomotion speed is lowest during the midstance phase and highest during the double-support phase.
This means that humans produce efficient walking through a pendular exchange of potential and kinetic energy while conserving mechanical energy~\cite{Cavagna1963, Cavagna1966,Cavagna1977}.
This is called the inverted pendulum mechanism~\cite{Ogihara2011}, and inverted pendulums have been widely used as the simplest model for the movement of the COM, when investigating the underlying mechanism of human walking~\cite{Alexander1980, Mochon1980a, Mochon1980b, Kuo2001, Kuo2001b, Srinivasan2007, Macdonald2014, Fujiki20150542}.

Passive dynamic walking, a popular dynamic system that is based on the inverted pendulum mechanism, was proposed by McGeer~\cite{McGeer1990, McGeer1993}.
This system walks down a shallow slope without an actuator or controller; it does this by balancing the energy dissipation due to foot contact with the energy generation due to the gravitational potential energy.
This walking behavior has various similarities to that of humans, and thus it has been a useful tool for elucidating the body's mechanical functions that produce energy-efficient walking~\cite{Asano2005, Bruijn2011, Chyou2011, Russell2005, Coleman1998, Collins2005, Collins2001, Goswami1998, Johnston2012, Kuo1999, Kurz2008, Kwan2007, Roos2010, Su2007}.

Due to the properties of the saddle point in the governing dynamics, an inverted pendulum falls down easily.
Therefore, a crucial issue is to clarify a stability mechanism for passive dynamic walking.
Garcia et al.~\cite{Garcia1998} used a simple compass-type model that incorporated the swing leg into the inverted pendulum (the simplest walking model), and they used a perturbation method to elucidate the generation of a stable limit cycle and the linear stability of the movement.
The stability characteristics of dynamical systems are determined by the basin of attraction of their attractors, as well as by their linear stability.
Schwab and Wisse~\cite{Schwab2001} investigated the basin of attraction of the simplest walking model and showed that it is very small and thin, and that it has a fractal-like shape; this explains the difficulty in producing stable passive dynamic walking.
There have been detailed studies of the stability and bifurcation of the simplest walking model and similar compass-type models~\cite{deBoer2010, Gritli2012a, Gritli2012b, LiYang2012}.
However, it remains unclear what mechanisms induce these geometric characteristics in the basin of attraction.

In the present study, we aim to clarify the mechanism that determines the geometric characteristics of the basin of attraction of the simplest walking model by considering the theory of dynamical systems and focusing on the saddle point that is inherent in the governing dynamics.
Because saddle points are embedded in general locomotor systems (they are not limited to passive dynamic walking), our results might contribute not only to elucidating the stability mechanism in passive dynamic walking, but also to improving the understanding of the stability mechanism in human walking and thus to producing design principles for the control of walking-support systems and biped robots.

\section{Methods}

\subsection{Model}

\begin{figure}[tbp]
  \centering
  \includegraphics[width=0.50\textwidth]{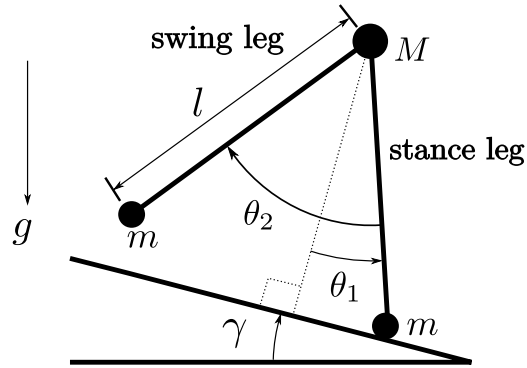}
  \caption{Simplest walking model}\label{fig:passive_config}
\end{figure}

In this study, we use the simplest walking model (Fig.~\ref{fig:passive_config}), introduced by Garcia et al.~\cite{Garcia1998}, for the dynamical analysis of passive dynamic walking.
This model has two legs (rigid links), each of length $l$, connected by a frictionless hip joint.
$\theta_1$ is the angle of the stance leg with respect to the slope normal, and $\theta_2$ is the angle between the stance leg and the swing leg.
The mass is located only at the hip and the feet; the hip mass is $M$, and the foot mass is $m$.
$g$ is the acceleration due to gravity.
This model walks on a slope of angle $\gamma$, without any control or input.
To simplify the analysis, we consider the limit case $\beta = m/M \to 0$, as in~\cite{Garcia1998}.

In this paper, we present a brief description of the model; please see~\cite{Garcia1998} for more details.

\subsubsection{Equations of motion for the swing phase}

The configuration of the model is described by two variables $(\theta_1, \theta_2)$, as indicated in Fig.~\ref{fig:passive_config}.
The equations of motion are as follows:
\begin{align}
  \ddot{\theta}_1 - \sin(\theta_1 - \gamma) &= 0, \label{eq:stance} \\
  (\cos\theta_2-1)\ddot{\theta}_1 + \ddot{\theta}_2 
  -\dot{\theta}_1^2\sin\theta_2 + \sin(\theta_2 - \theta_1 + \gamma)
                                            &= 0. \label{eq:swing}
\end{align}
Note that in these equations, we have already taken the limit as $\beta=m/M \to 0$, and we have nondimensionalized the equations using the time scale $\sqrt{l/g}$.

\subsubsection{Foot contact}

The swing foot contacts the slope when the following conditions are satisfied:
\begin{align}
  2 \theta_1 - \theta_2 &= 0,   \label{eq:cond1} \\
  \theta_1 &< 0,  \label{eq:cond2} \\
  2\dot{\theta}_1 - \dot{\theta}_2 &< 0. \label{eq:cond3} 
\end{align}
Conditions (\ref{eq:cond2}) and  (\ref{eq:cond3}) are used to ignore the foot scuffing when the swing leg moves forward.

We assume that foot contact is a fully inelastic collision (no slip, no bound) and that the stance foot lifts off the slope as soon as the swing foot hits the slope.
The relationship between the state just before foot contact $(\theta_1^-, \dot{\theta}_1^-, \theta_2^-, \dot{\theta}_2^-)$ and the state just after foot contact $(\theta_1^+, \dot{\theta}_1^+, \theta_2^+, \dot{\theta}_2^+)$ is as follows:
\begin{align}
  \begin{bmatrix}
    \theta_1^+ \\ \dot{\theta}_1^+ \\ \theta_2^+ \\ \dot{\theta}_2^+
  \end{bmatrix}
  =
  \begin{bmatrix}
    -\theta_1^- \\ \dot{\theta}_1^- \cos 2\theta_1^- \\ -2\theta_1^- \\ 
    \cos2\theta_1^-(1-\cos 2\theta_1^-)\dot{\theta}_1^-
  \end{bmatrix}. \label{eq:jump}
\end{align}
Since the state just after foot contact depends only on $(\theta_1^-, \dot{\theta}_1^-)$ and is independent of $(\theta_2^-, \dot{\theta}_2^-)$, it forms a two-dimensional surface in the four-dimensional phase space $(\theta_1, \dot{\theta}_1, \theta_2, \dot{\theta}_2)$.

\subsection{Structure of phase space by hybrid dynamics}

\begin{figure}[tbp]
  \centering
  \includegraphics[width=1.0\textwidth]{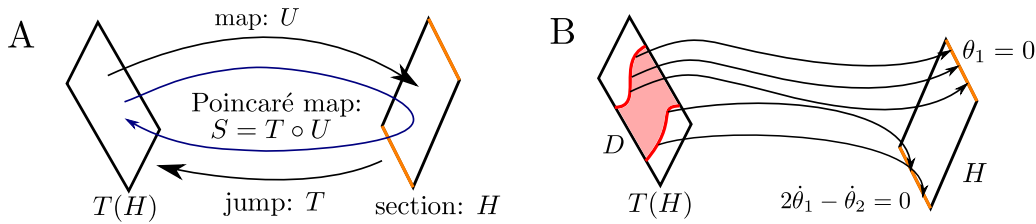} \\
  \caption{Structure of phase space $(\theta_1, \dot{\theta}_1, \theta_2, \dot{\theta}_2)$. A:\ Foot contact condition (section) $H$ bounded by two conditions (orange lines), the jump by foot contact $T$, the state just after the foot contact event $T(H)$, the map from $T(H)$ to $H$ by the equations of motion for the swing phase $U$, and the Poincar\'{e} map $S$ defined by $S = T \circ U$ on the Poincar\'{e} section $T(H)$. B:\ Domain $D$ (red region) bounded by the backward orbits of two boundaries of $H$ by the equations of motion for the swing phase (red lines)}\label{fig:domain}
\end{figure}

The simplest walking model is a hybrid system composed of the continuous dynamics during the swing phase and the discontinuous dynamics at foot contact.
This hybrid dynamic system determines the structure of the phase space, as shown in Fig.~\ref{fig:domain}A.
$H$ is the section of foot contact defined by the conditions (\ref{eq:cond1}), (\ref{eq:cond2}), and (\ref{eq:cond3}).
$T$ is the jump in the phase space from the state just before foot contact to the state just after foot contact, defined by the relationship (\ref{eq:jump}).
Therefore, the image of $T$, $T(H)$, is the region representing all states just after foot contact and a new step starts from $T(H)$.
$U$ is the map from the start of a step to the foot contact.
In other words, $U$ is the map from $T(H)$ to $H$, defined by the equations of motion (\ref{eq:stance}) and (\ref{eq:swing}).
The Poincar\'{e} map $S$ is defined by $S = T \circ U:T(H) \to T(H)$ on the Poincar\'{e} section $T(H)$.
This Poincar\'{e} map represents one step, and an attractor of the Poincar\'{e} map represents stable walking.
The basin of attraction of $S$ is the main topic of this paper.
$S$ is parameterized by one parameter $\gamma$, and Garcia et al.~\cite{Garcia1998} found that $S$ has an attracting fixed point at $0 < \gamma < 0.015$, and there is a period-doubling cascade to chaos for $0.015 < \gamma < 0.019$.

To investigate the basin of attraction, the domain of $T$ is important.
The map $S$ is not defined for all $T(H)$, since some initial conditions may cause the model to fall down.
We define the domain $D$ as the collection of initial conditions for which the model takes at least one step.
$D$ is in $T(H)$ and bounded, as shown in Fig.~\ref{fig:domain}B.
$H$ has two boundaries (orange lines) defined by $\theta_1=0$ and $2\dot{\theta}_1 - \dot{\theta}_2=0$ from the conditions (\ref{eq:cond2}) and (\ref{eq:cond3}), and the backward flows of these boundaries by the equations of motion (\ref{eq:stance}) and (\ref{eq:swing}) determine the boundaries of $D$ (red lines).

We also consider the sequence of inverse images of $D$, $S^{-n}(D)\ (n=1,2,\ldots)$.
These regions indicate the collections of initial conditions for which the model takes at least $(n+1)$ steps.
This sequence approximates the basin of attraction, and we investigate the mechanism by which the shape of the basin of attraction is formed from the geometric structure of these inverse images.

\subsection{Hyperbolicity and manifolds of the governed equations}

\begin{figure}[tbp]
  \centering
  \includegraphics[width=0.6\textwidth]{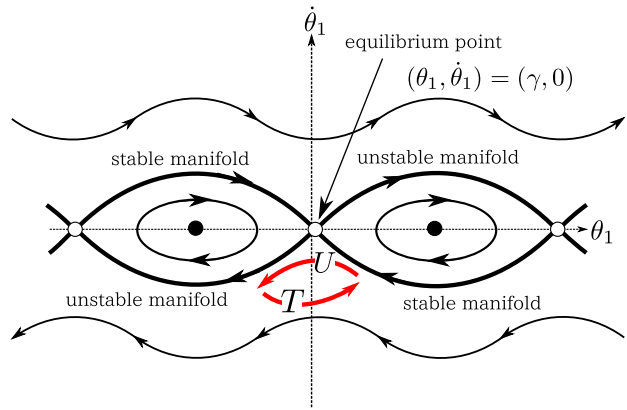}
  \caption{Phase diagram $(\theta_1, \dot{\theta}_1)$. The equilibrium point $(\theta_1, \dot{\theta}_1) = (\gamma, 0)$ is a saddle. Stable walking (red arrows) is obtained by the map $U$ and the jump $T$.
  }\label{fig:phase_space}
\end{figure}

The equations of motion (\ref{eq:stance}) and (\ref{eq:swing}) have an equilibrium point $(\theta_1, \dot{\theta}_1, \theta_2, \dot{\theta}_2) = (\gamma, 0, 0, 0)$.
At the equilibrium point, the legs remain upright.
The equilibrium point is deeply related to the geometric structure of the basin of attraction, which we will discuss in the following sections.
The eigenvalues of the linearized equations of motion at the equilibrium point are $\pm 1$ and $\pm i$, and the equilibrium point is a saddle-center with one stable direction, one unstable direction, and two neutral directions.

\begin{figure}[tpb]
  \centering
  \includegraphics[width=0.50\textwidth]{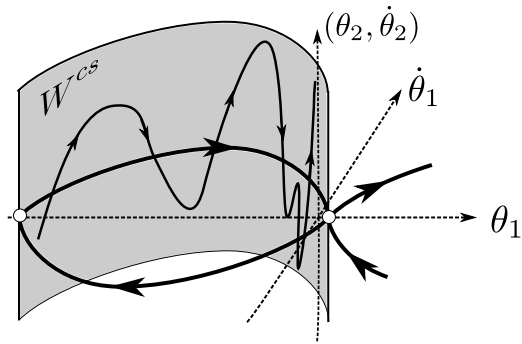}
  \caption{Center stable manifold $W^{cs}$}\label{fig:stable_mfd}
\end{figure}

The changes in the angle of the stance leg $\theta_1$ are governed by equation (\ref{eq:stance}) and are not affected by the movement of the swing leg $\theta_2$ (this is because we are considering the limiting case, $\beta \to 0$).
This equation for $\theta_1$ has a saddle equilibrium point at $(\theta_1, \dot{\theta}_1) = (\gamma, 0)$, as shown in Fig.~\ref{fig:phase_space}, similar to that of a single inverted pendulum.
In the phase diagram of $(\theta_1, \dot{\theta}_1)$ in Fig.~\ref{fig:phase_space}, bold lines going into $(\gamma, 0)$ are the stable manifold of the equilibrium point $W^s$, and the bold lines going out of the equilibrium point are the unstable manifold of the equilibrium point $W^u$.
In the phase space of four variables, $(\theta_1, \dot{\theta}_1, \theta_2, \dot{\theta}_2)$, $W^s\times\Real^2$ and $W^u\times\Real^2$ are the center stable manifold and the center unstable manifold, respectively, and we denote them by $W^{cs}$ and $W^{cu}$.
An orbit on $W^{cs}$ behaves as follows (Fig.~\ref{fig:stable_mfd}):
\begin{itemize}
\item An orbit starting from a point on $W^{cs}$ never goes outside of $W^{cs}$;
\item An orbit starting from a point on $W^{cs}$ converges to $(\gamma, 0)\times \Real^2$ as $\mbox{(time)} \to +\infty$.
\end{itemize}
An orbit on $W^{cu}$ behaves in the same way as one on $W^{cs}$, as $\mbox{(time)} \to -\infty$.

\subsection{Computation of regions and manifolds}

To investigate the geometric structure of the basin of attraction, we compute the domain $D$, the sequence of inverse images of the domain $S^{-n}(D)\ (n=1,2,\ldots)$, the basin of attraction $B$, and the intersection of $W^{cs}$ and $T(H)$.
All of these sets are defined on $T(H)$.
The jump map (\ref{eq:jump}) shows that $T(H)$ is a two-dimensional surface in the four-dimensional phase space $\Real^4$, and each point on $T(H)$ is uniquely determined by two variables $(\theta_1, \dot{\theta}_1)$.
Therefore, we use the coordinates $(\theta_1, \dot{\theta}_1)$ to describe this region.
Because $D$ is a set of initial points on $T(H)$ that reach $H$ through the equations of motion (\ref{eq:stance}) and (\ref{eq:swing}) (Fig.~\ref{fig:domain}), we can compute $D$ by numerically integrating the equations.
We can also compute $S^{-n}(D)\ (n=1,2,\ldots)$ \ in a similar way.
We approximate $B$ as a set of initial points on $T(H)$ that allows the model to take a sufficient number of steps.
More specifically, we compared the results of 50 and 200 steps and used the initial points when the two results were identical.
We can also compute the intersection of $W^{cs}$ and $T(H)$ from the fact that $W^{cs}$ is a separatrix, as shown in Fig.~\ref{fig:phase_space}.

\section{Results} 

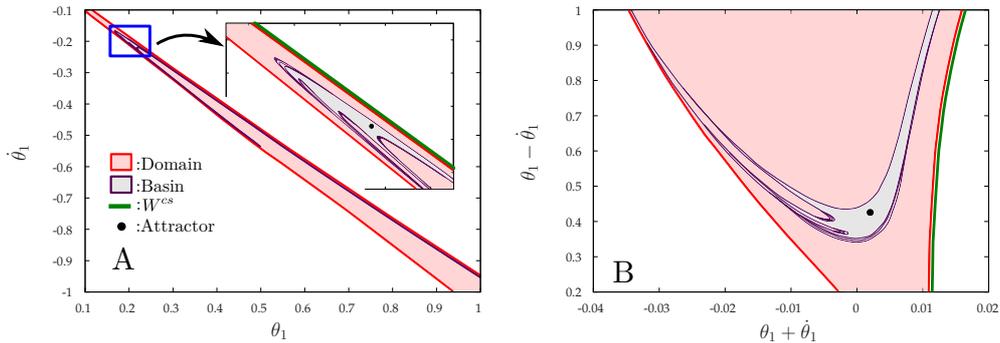
\begin{figure}[t]
  \centering 
  \resizebox{0.49\textwidth}{!}{\input{basin_domain_1z.tex}}
  \resizebox{0.49\textwidth}{!}{\input{basin_domain_3x.tex}} 
  
  \caption{Geometric characteristics of the basin of attraction for $\gamma=0.011$. A:\ Domain $D$ and basin of attraction $B$ on $(\theta_1, \dot{\theta}_1)$. The region enclosed by the blue box is magnified. B:\ Rotated view using $\theta_1 + \dot{\theta}_1$ and $\theta_1 - \dot{\theta}_1$ for the axes.}\label{fig:domain_basin_0_011}
\end{figure}

\subsection{Geometric characteristics of the domain and basin of attraction}

Figure~\ref{fig:domain_basin_0_011}A shows the domain $D$ and the basin of attraction $B$ on $T(H)$ with $\gamma=0.011$.
Both $D$ and $B$ are very thin in the space of $(\theta_1, \dot{\theta}_1)$.
To clearly see the geometrical details, we rotated the figure using $\theta_1 + \dot{\theta}_1$ and $\theta_1 - \dot{\theta}_1$ for the axis in Fig.~\ref{fig:domain_basin_0_011}B.
The intersection of the center-stable manifold $W^{cs}$ and $T(H)$ is shown by a green line in Figs.~\ref{fig:domain_basin_0_011}A and B.
From these figures, we found that $D$ had the following properties:
\begin{itemize}
\item $D$ has a long, thin region in $(\theta_1, \dot{\theta}_1)$;
\item Two boundaries of $D$ are almost parallel, and one of them is very close to $W^{cs}$.
\end{itemize}
We also found the following properties for $B$:
\begin{itemize}
\item $B$ is located inside $D$ and is thinner than $D$;
\item $B$ is V-shaped;
\item There are fractal-like slits in $B$ and a stripe pattern in the cusp of the V-shaped region.
\end{itemize}

\subsection{Geometric characteristics of the inverse images of the domain}

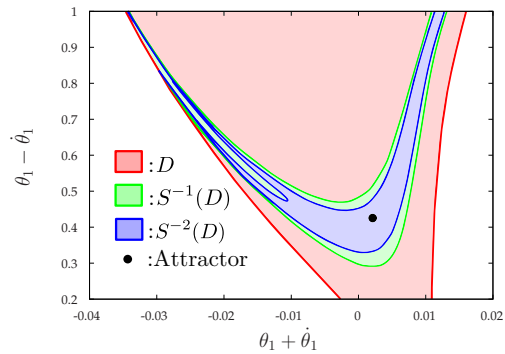
\begin{figure}[t]
  \centering
  \resizebox{0.50\textwidth}{!}{\input{domain123.tex}}

  \caption{Domain $D$ and inverse images $S^{-1}(D)$ and $S^{-2}(D)$}\label{fig:domain123}
\end{figure}

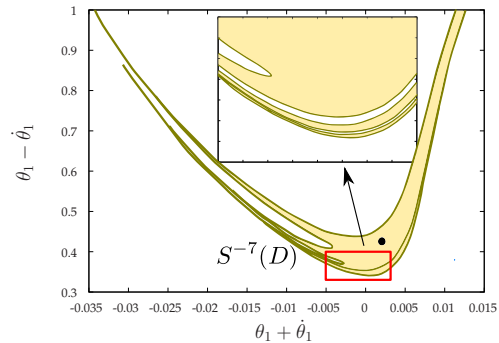
\begin{figure}[t]
  \centering
  \resizebox{0.49\textwidth}{!}{\input{domain8_8L.tex}}
  \caption{Inverse image $S^{-7}(D)$. The region enclosed by the red box is magnified.}\label{fig:domain8}
\end{figure}

To investigate how to generate these geometric characteristics of $B$ from $D$, we calculated the inverse images of $D$, $S^{-n}(D) \ (n=1,2,\ldots)$.
Figure~\ref{fig:domain123} shows $D$, $S^{-1}(D)$, and $S^{-2}(D)$.
We found the following:
\begin{itemize}
\item $S^{-1}(D)$ is contained in $D$ and is V-shaped;
\item $S^{-2}(D)$ is located inside $S^{-1}(D)$, is V-shaped, and has a slit.
\end{itemize}
Figure~\ref{fig:domain8} shows the inverse image $S^{-7}(D)$.
A stripe pattern appears through the sequence of $S^{-n}(D)\ (n=1,2,\ldots)$.

As shown in the above figures, $B$ is constructed by the sequence of the backward images $S^{-n}(D)\ (n=1,2,\ldots)$.
Therefore, we can find the construction mechanism of the shape of the basin of attraction from the backward images; we will discuss this in Section~\ref{sec:discussion}.

\subsection{Dependence of the slope angle}

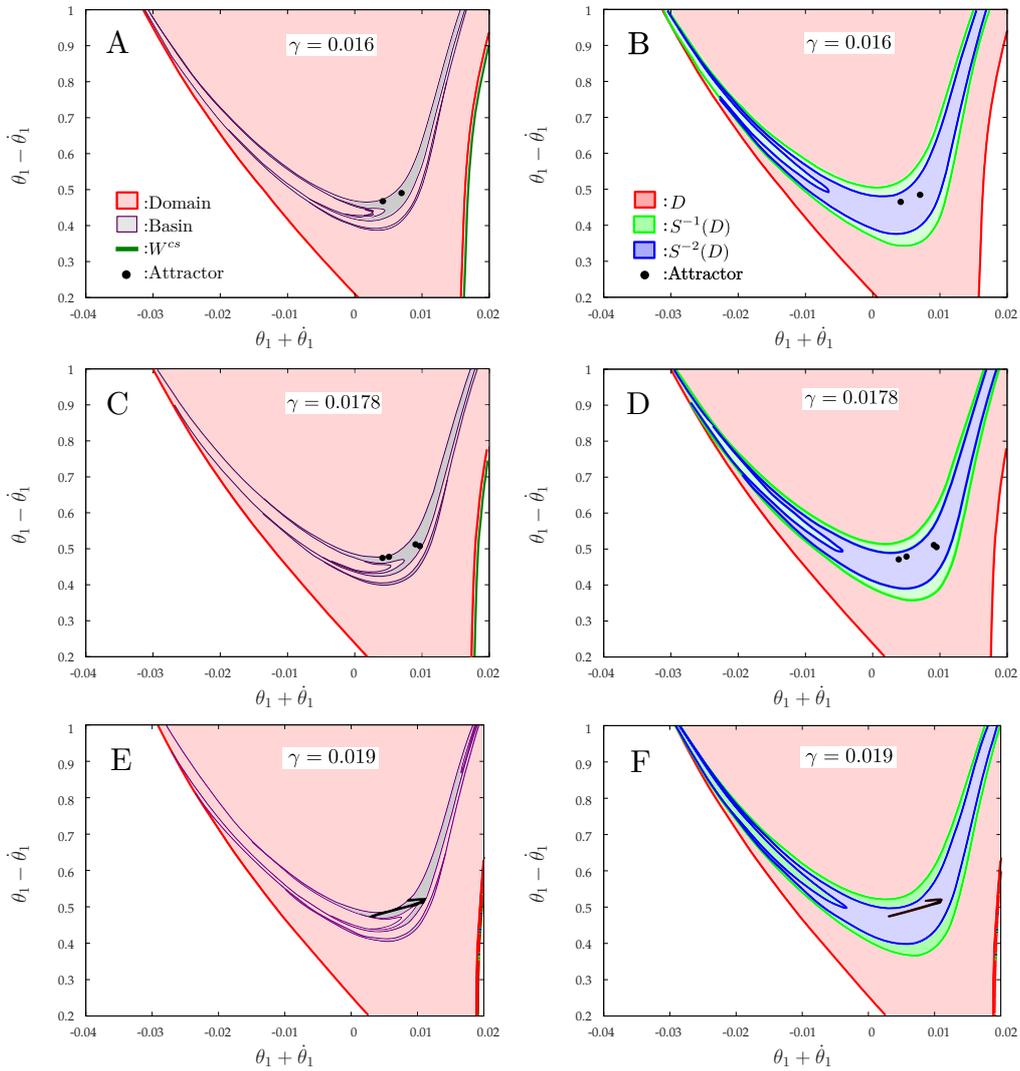
\begin{figure*}
  \resizebox{0.5\textwidth}{!}{\input{0_0160_basin.tex}}
  \resizebox{0.5\textwidth}{!}{\input{0_0160.tex}}
  \resizebox{0.5\textwidth}{!}{\input{0_01780_basin.tex}} 
  \resizebox{0.5\textwidth}{!}{\input{0_01780.tex}}
  \resizebox{0.5\textwidth}{!}{\input{0_01903_basin.tex}}
  \resizebox{0.5\textwidth}{!}{\input{0_01903.tex}}
  \caption{Dependence of the parametric characteristics on the parameter $\gamma$. $\gamma=0.016$ for A and B, 0.0178 for C and D, and 0.019 for E and F. A, C, and E show the domain $D$, the basin of attraction $B$, and the intersection with $W^{cs}$. B, D, and F show the domain $D$, and the inverse images $S^{-1}(D)$ and $S^{-2}(D)$.}\label{fig:params}
\end{figure*}

To examine the dependence of the geometric characteristics on the parameter $\gamma$, we calculated the domain $D$, the basin of attraction $B$, and the inverse images $S^{-1}(D)$ and $S^{-2}(D)$ for various values of $\gamma$.
Figures~\ref{fig:params}A and B show the results for $\gamma = 0.016$ and have two attracting points, which correspond to a stable period-2 gait.
Figures~\ref{fig:params}C and D show the results for $\gamma=0.0178$ and have four attracting points, which correspond to a stable period-4 gait.
Figures~\ref{fig:params}E and F show the results for $\gamma=0.019$ and have a chaotic attractor.
Although the characteristics of the attractor change through period-doubling bifurcations depending on $\gamma$, the shapes of $B$, $D$, $S^{-1}(D)$, and $S^{-2}(D)$ do not change very much.
When $\gamma$ is larger than 0.019, the attractor disappears.
We will discuss this ``attractor crisis'' phenomenon in Section~\ref{sec:discussion}\ref{subsec:crisis}.

\section{Discussion}\label{sec:discussion}

\begin{figure}[tbp]
  \centering
  \begin{minipage}{0.45\textwidth}
    \centering
    \includegraphics[width=\hsize]{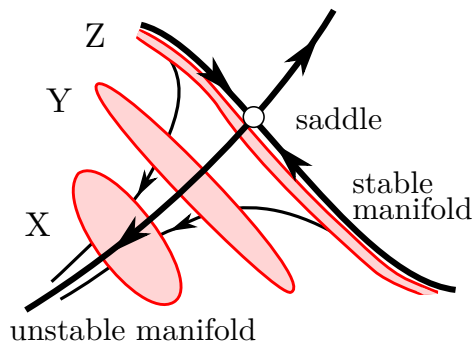}
  \end{minipage} 
  
  \caption{$\lambda$-lemma. Region X moves and is deformed to thinner regions Y and Z by the backward flow. In our model, the foot contact section $H$ and domain $D$ correspond to the regions X and Z, respectively.
  }\label{fig:lambda_lemma}
\end{figure}

\subsection{Why the domain is thin along the center stable manifold}

The domain $D$ is very thin, as shown in Fig.~\ref{fig:domain_basin_0_011}A.
This is related to the $\lambda$-lemma, one of the most important theorems in the theory of dynamical systems~\cite{Robinson2008}.
From this theorem, we can say the following:
\begin{itemize}
\item A region intersecting the unstable manifold of a saddle equilibrium point moves toward the stable manifold of the saddle when the region is moved by the flow that is backward in time;
\item When the region comes close to the stable manifold, the region becomes thin due to the hyperbolic behavior near the saddle.
\end{itemize}
Figure~\ref{fig:lambda_lemma} illustrates how a region X moves and is deformed into thinner regions Y and Z through the backward flow.
As shown in Fig.~\ref{fig:domain}B, $D$ is obtained by the intersection of $T(H)$ and the backward orbit whose initial point is in $H$.
Therefore, $D$ becomes thin along the center stable manifold, as shown in Figs.~\ref{fig:domain_basin_0_011}A and \@\ref{fig:lambda_lemma}.

\subsection{Why the inverse image of the domain is V-shaped}

Since the sequence of the inverse images $S^{-n}(D)\ (n=1,2,\ldots)$ approximates the basin of attraction, it is important to clarify how the geometric structure of the inverse images is constructed, and hence clarify the shape of the basin of attraction.

First, we discuss why $S^{-1}(D)$ is V-shaped in the thin region $D$, as shown in Fig.~\ref{fig:domain123}.
Since $S^{-1}(D) = U^{-1}(T^{-1}(D))$, we show below how the shape of $T^{-1}(D)$ is deformed by $U^{-1}$.

\subsubsection{Shape of $T^{-1}(D)$}

\begin{figure}[tbp]
  \centering
  \includegraphics[width=0.50\textwidth]{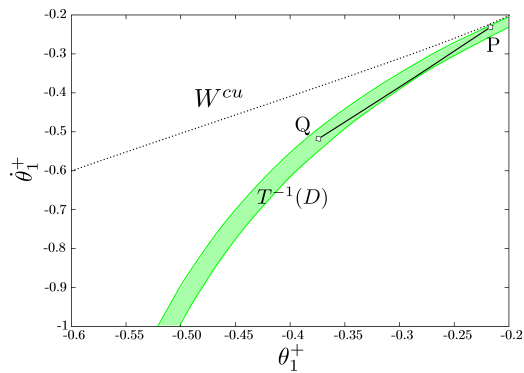} 
  \caption{$T^{-1}(D)$ in $(\theta_1^+, \dot{\theta}_1^+)$}\label{fig:T-1D}
\end{figure}

From (\ref{eq:jump}), $T^{-1}(D)$ is described as:
\begin{align}
  \{&(-\theta_1^+, -2\theta_1^+, \dot{\theta}_1^+  \sec 2\theta_1^+,
      \dot{\theta}_2^-) 
      \mid  \nonumber \\ 
    &(\theta_1^+, \theta_2^+, \dot{\theta}_1^+, \dot{\theta}_2^+) \in
      D , \dot{\theta}_2^- \in {\mathbb R} 
      \}. \label{eq:Tinv}
\end{align}
This is explained by $\theta_1^+$, $\dot{\theta}_1^+$, and $\dot{\theta}_2^-$.
Figure~\ref{fig:T-1D} shows numerically obtained results for $T^{-1}(D)$ in $(\theta_1^+, \dot{\theta}_1^+)$.
The region is thin and curved.
The relative positions of $T^{-1}(D)$ and $W^{cu}$ are important for a V-shaped $S^{-1}(D)$, as explained below.

\subsubsection{Deformation of $T^{-1}(D)$ by $U^{-1}$}

We denote the solution of equations of motion \eqref{eq:stance} and \eqref{eq:swing} by $\Theta(t) = (\theta_1(t), \theta_2(t), \dot{\theta}_1(t), \dot{\theta}_2(t))$.
From the definition of $U$, for a point $\Theta(0) \in T^{-1}(D) \subset H$, there exists $\Delta > 0$ such that
\begin{align}
  \Theta(-\Delta) = U^{-1}(\Theta(0)) \in T(H) \label{eq:U-no1}
\end{align}
holds (to analyze $U^{-1}$, it is natural to consider this in terms of negative time: $-\Delta$ (Fig.~\ref{fig:domain}A)).
$\Theta(0)$, $\Theta(-\Delta)$, and $\Delta$ correspond to the state just before foot contact, the state just after foot contact, and the duration of a step, respectively.
$\Theta(-\Delta)$ gives the deformation of $T^{-1}(D)$ by $U^{-1}$.

Since $\Theta(-\Delta)$ is contained in $T(H)$, the following equations hold:
\begin{align}
  2\theta_1(-\Delta) &= \theta_2(-\Delta), \label{eq:start_1} \\
  \dot{\theta}_2(-\Delta) &= \dot{\theta}_1(-\Delta)(1 -
                            \cos2\theta_1(-\Delta)), \label{eq:start_2} \\
  \theta_1(-\Delta) &> 0. \label{eq:start_3}
\end{align}
The equality (\ref{eq:start_2}) comes from $\dot{\theta}_2^+ = \dot{\theta}_1^- \cos2\theta_1^-(1-\cos 2\theta_1^-)$ and $\dot{\theta}_1^+ = \dot{\theta}_1^- \cos 2\theta_1^-$ in (\ref{eq:jump}).
In addition, since $\Theta(-\Delta)$ is contained in $H$, the following equation holds:
\begin{align}
  2\theta_1(0) &= \theta_2(0). \label{eq:end_1}
\end{align}

To approximately solve \eqref{eq:U-no1}, we linearize equations (\ref{eq:stance}) and (\ref{eq:swing}) around $(\gamma, 0,0,0)$ by using
\begin{align*}
  \ddot{\theta}_1 &= \theta_1 - \gamma, \\
  \ddot{\theta}_2 &= -(\theta_2 - \theta_1 + \gamma).
\end{align*}
The solution is:
\begin{align}
  \theta_1 &= \gamma + C_1 \exp(t) + C_2 \exp(-t), \label{eq:approx_saddle} \\
  \theta_2 - (\theta_1 - \gamma)/2 &= K \cos(t + \phi), \label{eq:approx_pendulum}
\end{align}
where $C_1, C_2, K$, and $\phi$ are the integration constants ($0 \leq \phi < 2\pi$).
This shows that the motion for the swing phase consists of two dynamic components: an inverted pendulum (\ref{eq:approx_saddle}) and a normal pendulum~\eqref{eq:approx_pendulum}.
$C_1$ and $C_2$ are determined by the initial conditions of $\theta_1$ and $\dot{\theta}_1$, as follows:
\begin{equation}
  \begin{aligned}
    C_1 &= (\theta_1(0) - \gamma + \dot{\theta}_1(0))/2, \\
    C_2 &= (\theta_1(0) - \gamma - \dot{\theta}_1(0))/2.
  \end{aligned}
  \label{eq:c1_c2}
\end{equation}
In contrast, $K$ and $\phi$ are determined by the initial conditions of $\theta_1, \theta_2, \dot{\theta}_1$, and $\dot{\theta}_2$.
In the linearized equations, $W^{cs}$ and $W^{cu}$ are approximated by $E^{cs} = \{(\theta_1, \theta_2, \dot{\theta}_1, \dot{\theta}_2) \mid \theta_1 - \gamma = -\dot{\theta}_1 \}$ and $E^{cu} = \{(\theta_1, \theta_2, \dot{\theta}_1, \dot{\theta}_2) \mid \theta_1 - \gamma = \dot{\theta}_1 \}$, respectively.

From \eqref{eq:start_1}, \eqref{eq:start_2}, \eqref{eq:end_1}, \eqref{eq:approx_saddle}, and \eqref{eq:approx_pendulum}, we have the following system of equations:
\begin{align}
  &\theta_1(-\Delta) = C_1 \exp(-\Delta) + C_2 \exp\Delta + \gamma,
    \label{eq:theta1_D} \\
  &\dot{\theta}_1(-\Delta) = C_1 \exp(-\Delta) - C_2 \exp\Delta,
    \label{eq:dtheta1_D} \\
  &K\cos(-\Delta + \phi) = 
    3\theta_1(-\Delta)/2  + \gamma/2, \label{eq:KcosDelta} \\
  &K\sin(-\Delta + \phi) = 
    -\dot{\theta}_1(-\Delta)(1/2 - \cos2\theta_1(-\Delta)), \label{eq:KsinDelta} \\
  &K\cos\phi = 3/2 \cdot (C_1 + C_2 + 4\gamma/3), \label{eq:Kcos0}
\end{align}
where $\Delta, \phi, K, \theta_1(-\Delta), \mbox{ and } \dot{\theta}_1(-\Delta)$ are unknown variables ($C_1$ and $C_2$ are determined in \eqref{eq:c1_c2} from $(\theta_1(0), \dot{\theta}_1(0))$).
We can compute $\Theta(-\Delta)$ \eqref{eq:U-no1} from $(\theta_1(0), \dot{\theta}(0))$ by solving (\ref{eq:theta1_D}-\ref{eq:Kcos0}).

\begin{figure}[t]
  \centering
  \includegraphics[width=0.60\textwidth]{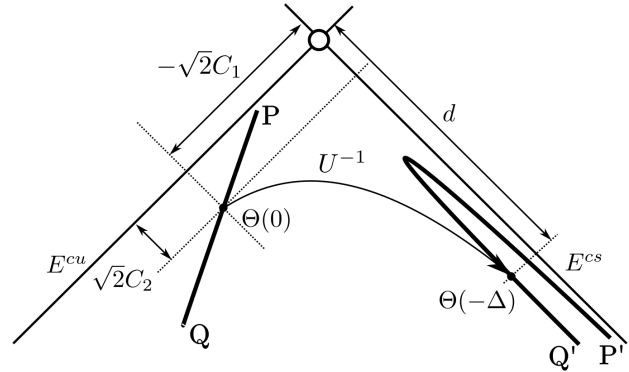}
  \caption{Translation of $\Theta(0)$ to $\Theta(-\Delta)$ by $U^{-1}$ for the linearized equations of motion}\label{fig:c1_c2_d_U}
\end{figure}

\begin{figure}[t]
  \centering
  \includegraphics[width=0.70\textwidth]{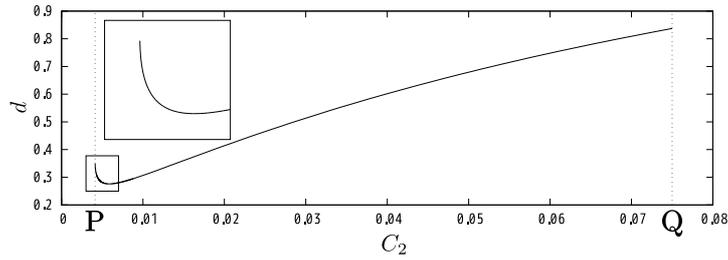}
  \caption{$d$ versus $C_2$}\label{fig:c2_d}
\end{figure}

To show how $U^{-1}$ deforms $T^{-1}(D)$, we use the approximated solution to show the relationship between the distance between $T^{-1}(D)$ and $W^{cu}$ and that between $U^{-1}(T^{-1}(D))$ and $W^{cu}$.
We use a line segment PQ ($\dot{\theta}_1 = 0.86 \theta_1 - 0.037, -0.364 < \theta_1 < -0.199 $) within $T^{-1}(D)$, as shown in Fig.~\ref{fig:T-1D}.
This is moved to the curve ${\rm P'Q'} = U^{-1}({\rm PQ})$, which approximates $U^{-1}(T^{-1}(D))$, as shown in Fig.~\ref{fig:c1_c2_d_U}.
When we take $\Theta(0)$ on PQ, we obtain the relation between $C_1$ and $C_2$ from (\ref{eq:c1_c2}), as follows:
\begin{equation}
  \begin{aligned}
    &C_1 = -0.2 + 13.3C_2,  \\
    &0.0038 < C_2 < 0.075 \ (-0.45 < C_1 < -0.22).
  \end{aligned}
  \label{eq:c1_c2_onPQ}
\end{equation}
$\sqrt{2}C_2$ corresponds to the distance between $\Theta(0)$ and $E^{cu}$, and the point on PQ is parameterized by $C_2$.
Figure~\ref{fig:c2_d} shows a graph of $d$ versus $C_2$, where $d$ is the distance between $\Theta(-\Delta)$ and $E^{cu}$ (Fig.~\ref{fig:c1_c2_d_U}).
In this figure, $C_2$ does not begin at 0, but at around 0.004; this is due to a singularity, which will be discussed in the next section.
The graph of $d$ is convex, which means that $T^{-1}(D)$ is nonuniformly deformed by $U^{-1}$.
Both the part of $T^{-1}(D)$ very close to $E^{cu}$ (around P) and the part of $T^{-1}(D)$ far from $E^{cu}$ (around Q) are strongly deformed by $U^{-1}$, and the part in between is relatively weakly deformed.
Therefore, $U^{-1}(T^{-1}(D)) = S^{-1}(D)$ is V-shaped, as shown in Figs.~\ref{fig:domain_construction_1}A and B.
Below, we will present a mathematical analysis to explain why $d$ is convex.

\subsubsection{Mathematical analysis of the deformation by $U^{-1}$}

\begin{figure}[t]
  \centering
  \includegraphics[width=0.7\textwidth]{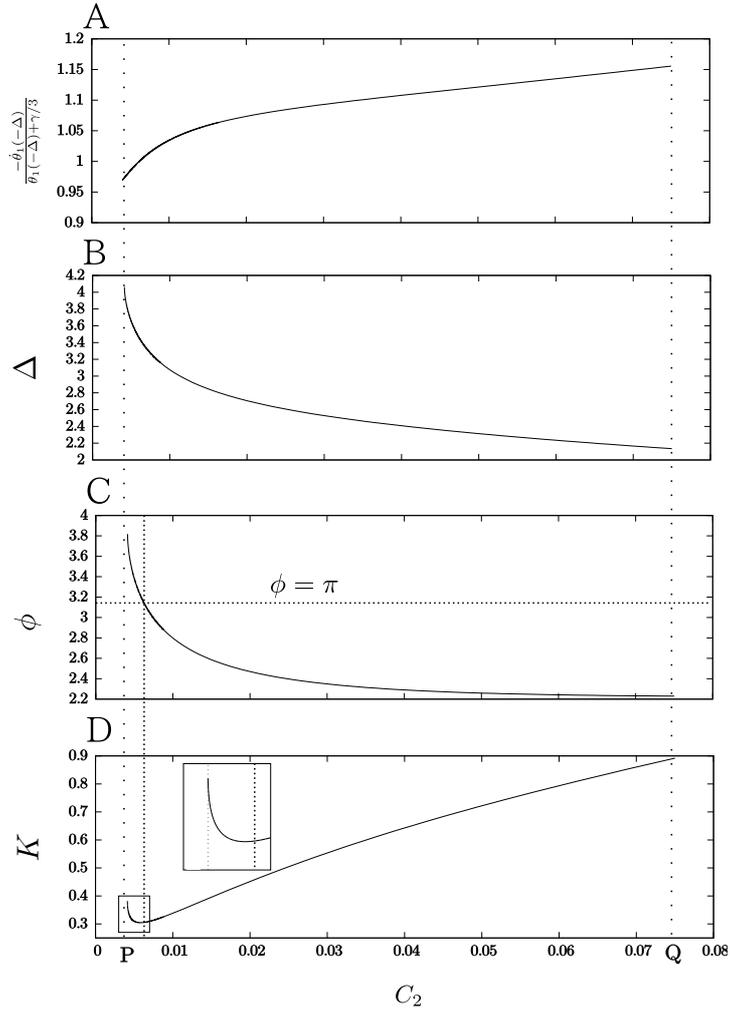}
  \caption{$\frac{-\dot{\theta}_1(-\Delta)}{\theta_1(-\Delta) + \gamma/3}$ (A), $\Delta$ (B), $\phi$ (C), and $K$ (D) versus $C_2$. In (D), $K$ is convex, with a minimum when $\phi\approx\pi$.}\label{fig:vs_c2}
\end{figure}

First, we consider the case where $\Theta(0)$ is close to $E^{cu}$.
In this case, $C_2$ is small.
From (\ref{eq:KcosDelta}) and (\ref{eq:KsinDelta}), we obtain
\begin{align}
  \tan(-\Delta + \phi) = 
  \frac{-\dot{\theta}_1(-\Delta)\cdot (1 - 2\cos 2\theta_1(-\Delta))}{3\theta_1(-\Delta) + \gamma}.\label{eq:tanNo1}
\end{align}
Since $\Theta(-\Delta) \in S^{-1}(D) \subset D$ and $D$ is very thin along $W^{cs} \approx E^{cs}$, we obtain
\begin{align*}
  \theta_1(-\Delta) \approx -\dot{\theta}_1(-\Delta) + \gamma
\end{align*}
and
\begin{align}
  \frac{-\dot{\theta}_1(-\Delta)}{\theta_1(-\Delta) + \gamma/3}
  \approx 1 \label{eq:dt_D_div_dt}
\end{align}
when $|\theta_1(-\Delta)|$ is sufficiently larger than $\gamma$.
In fact, the left-hand side of (\ref{eq:dt_D_div_dt}) is 0.95--1.16 on PQ, as shown in Fig.~\ref{fig:vs_c2}A.

From (\ref{eq:tanNo1}) and (\ref{eq:dt_D_div_dt}), we approximately have
\begin{align}
  \tan(-\Delta + \phi) = 2(1/2 - \cos 2\theta_1(-\Delta))/3.
  \label{eq:tanDp}
\end{align}
Since $-1 \leq \cos 2\theta_1(-\Delta) \leq 1$, we have
\begin{align}
  -1 \leq -1/3 \leq \tan(-\Delta + \phi) \leq 1\label{eq:tanDp13}
\end{align}
and $|\cos(-\Delta + \phi)| \geq 1/\sqrt{2}$.
From (\ref{eq:start_3}) and (\ref{eq:KcosDelta}), we have $\cos(-\Delta + \phi) \geq 0$ and
\begin{align}
  \cos(-\Delta + \phi) \geq 1/\sqrt{2}. \label{eq:cosRange}
\end{align}
Therefore, from (\ref{eq:theta1_D}), (\ref{eq:KcosDelta}), (\ref{eq:Kcos0}), and $C_1 < 0$, we have
\begin{align*}
  C_2 \exp \Delta
  &= \theta_1(-\Delta) - C_1 \exp(-\Delta) - \gamma \\
  &\geq \theta_1(-\Delta) - \gamma \\
  &= (2/3)\cdot K \cos(-\Delta + \phi) - 4\gamma/3 \\
  & \geq (2/3)\cdot K \cdot(1/\sqrt{2})
    - 4\gamma/3 \\
  & = (\sqrt{2}/3)K - 4\gamma/3 \\
  & = 1/\sqrt{2} \cdot (C_1 + C_2 + 4\gamma/3)/ \cos \phi - 4\gamma/3 \\
  &\geq (1/\sqrt{2})\cdot|C_1 + C_2 + 4\gamma/3| - 4\gamma/3.
\end{align*}
Since from (\ref{eq:c1_c2_onPQ}) we have that $C_1$ is much larger than both $|\gamma|$ and $|C_2|$, $\Delta$ goes to $+\infty$ as $C_2 \to 0$.
Therefore, $\Delta$ increases as $C_2$ decreases, as shown in Fig.~\ref{fig:vs_c2}B.

From (\ref{eq:cosRange}), we obtain \[ -\pi/4 < -\Delta + \phi< \pi/4.\]
Since $\Delta - \pi/4 < \phi < \Delta + \pi/4$ and $\Delta \to \infty$ as $C_2 \to 0$, $\phi$ goes to $\infty$ as $C_2 \to 0$ and $\phi$ increases when $C_2$ decreases as shown Fig.~\ref{fig:vs_c2}C.
In addition, $\pi/2 < \phi < 3\pi/2$ since $K\cos\phi = 3/2 \cdot (C_1 + C_2 + 4\gamma/3) < 0$ from (\ref{eq:c1_c2_onPQ}), as shown Fig.~\ref{fig:vs_c2}C.
From the above two facts, the solution of (\ref{eq:theta1_D}-\ref{eq:Kcos0}) has a singularity when $C_2$ is small.

As shown in Fig.~\ref{fig:vs_c2}C, $\phi$ increases as $C_2$ decreases.
As $\phi$ monotonically increases through $\pi$, $\left|\cos \phi\right|$ monotonically increases when $\phi$ is less than $\pi$ and decreases when it is greater than $\pi$.
In addition, the changes in $C_1$ and $C_2$ are much smaller than that of $\cos \phi$.
This is why $K = 3/2 \cdot (C_1 + C_2 + 4\gamma/3) / \cos\phi$ is minimized at $\phi \approx \pi$, as shown Fig.~\ref{fig:vs_c2}D.
This also causes $\theta_1(-\Delta) = 2|K \cos(-\Delta + \phi)|/3 - \gamma/3$ and the curve of $d = \sqrt{2}\cdot(\theta_1(-\Delta))$ to be minimized near that point.
These factors combine to cause $d$ to be convex when $C_2$ is small.
This is the reason why $d$ has a convex for a small $C_2$.

We can explain this phenomenon intuitively from a dynamic viewpoint.
The gait is generated by the coordination between the inverted pendulum behavior of the stance leg \eqref{eq:approx_saddle} and the normal pendulum behavior of the swing leg \eqref{eq:approx_pendulum}.
When $C_2$ is very small, $(\theta_1(t), \dot{\theta}(t))$ passes through the small neighborhood of the saddle equilibrium point $(\gamma, 0)$.
Therefore the inverted pendulum behavior~\eqref{eq:approx_saddle} becomes very slow, and $\Delta$, which corresponds to step duration, becomes very large.
In contrast, the normal pendulum behavior~\eqref{eq:approx_pendulum} does not become slow, since the angular velocity is $1$ independent of $C_2$, as shown in~\eqref{eq:approx_pendulum}.
To deal with this, the initial phase $\phi$ becomes larger than $\pi$, there in an increase in the magnitude of the swing leg motion represented by $K$, and $d$ increases.

\begin{figure}[t]
  \centering
  \includegraphics[width=\hsize]{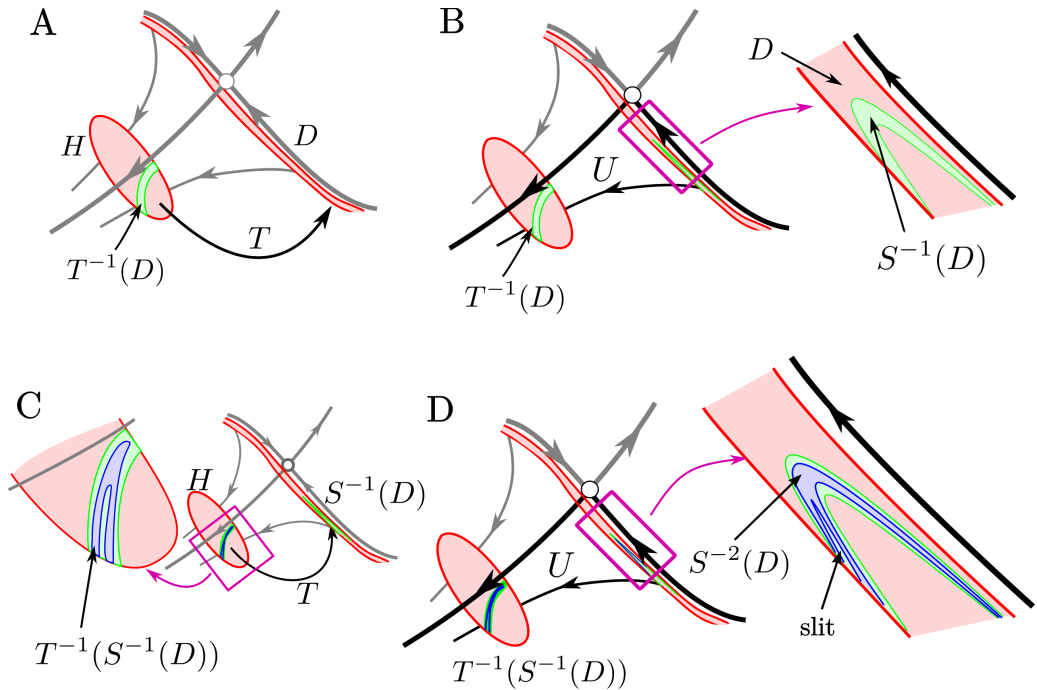}
  \caption{
    Formation of the V-shape and split in the inverse images of the domain. A:~$T^{-1}(D)$ is obtained by the inverse image of $D$ and becomes a thin and curved region in $H$. B:~$T^{-1}(D)$ is moved and deformed by the backward flow to $S^{-1}(D)=U^{-1}(T^{-1}(D))$. C:~$T^{-1}(S^{-1}(D))$ is obtained by the inverse image of $S^{-1}(D)$. D:~$T^{-1}(S^{-1}(D))$ is moved and deformed by the backward flow to $S^{-2}(D)=U^{-1}(T^{-1}(S^{-1}(D)))$.
  }\label{fig:domain_construction_1}
\end{figure}

Next, we consider the case of large $C_2$.
In this case, $\Delta$ is not very large, since the orbit does not pass through the neighborhood of the saddle equilibrium point.
Therefore, $C_2$ dominates $d$, since $d/\sqrt{2} + \gamma = \theta_1(-\Delta) = C_1 \exp(-\Delta) + C_2 \exp\Delta +\gamma$, and so $d$ increases as $C_2$ increases.

If we integrate the above two cases, the part of $T^{-1}(D)$ very close to $E^{cu}$ (around P) and the part of $T^{-1}(D)$ far from $E^{cu}$ (around Q) are both strongly deformed, and the region between them is weakly deformed.
This is because $S^{-1}(D)$ is V-shaped.
This mechanism is illustrated in Figures~\ref{fig:domain_construction_1}A and B.

\subsection{Why the inverse images have slits and stripe patterns}

The inverse image $S^{-2}(D)$ has a slit due to a mechanism similar to the one discussed in the previous section.
Since $S^{-2}(D) = U^{-1}(T^{-1}(S^{-1}(D)))$, we consider two steps, $T^{-1}(S^{-1}(D))$ and $U^{-1}(T^{-1}(S^{-1}(D)))$.
Here, $T^{-1}(S^{-1}(D))$ is obtained by the backward image of $S^{-1}(D)$, and it is contained in $T^{-1}(D)$, as shown by the blue region in Fig.~\ref{fig:domain_construction_1}C.
This region is moved by $U^{-1}$; it is expanded along the direction of the stable manifold and contracted along the direction of the unstable manifold, as with $T^{-1}(D)$.
As a result, $S^{-2}(D)$ becomes V-shaped with a slit, as shown in Fig.~\ref{fig:domain_construction_1}D.
We can also give a similar explanation for the stripe pattern, which is formed by the repeated expansion of nested regions.

\subsection{Why the attractor disappears}\label{subsec:crisis}

\begin{figure}[!t]
  \centering
  \includegraphics[width=0.49\textwidth]{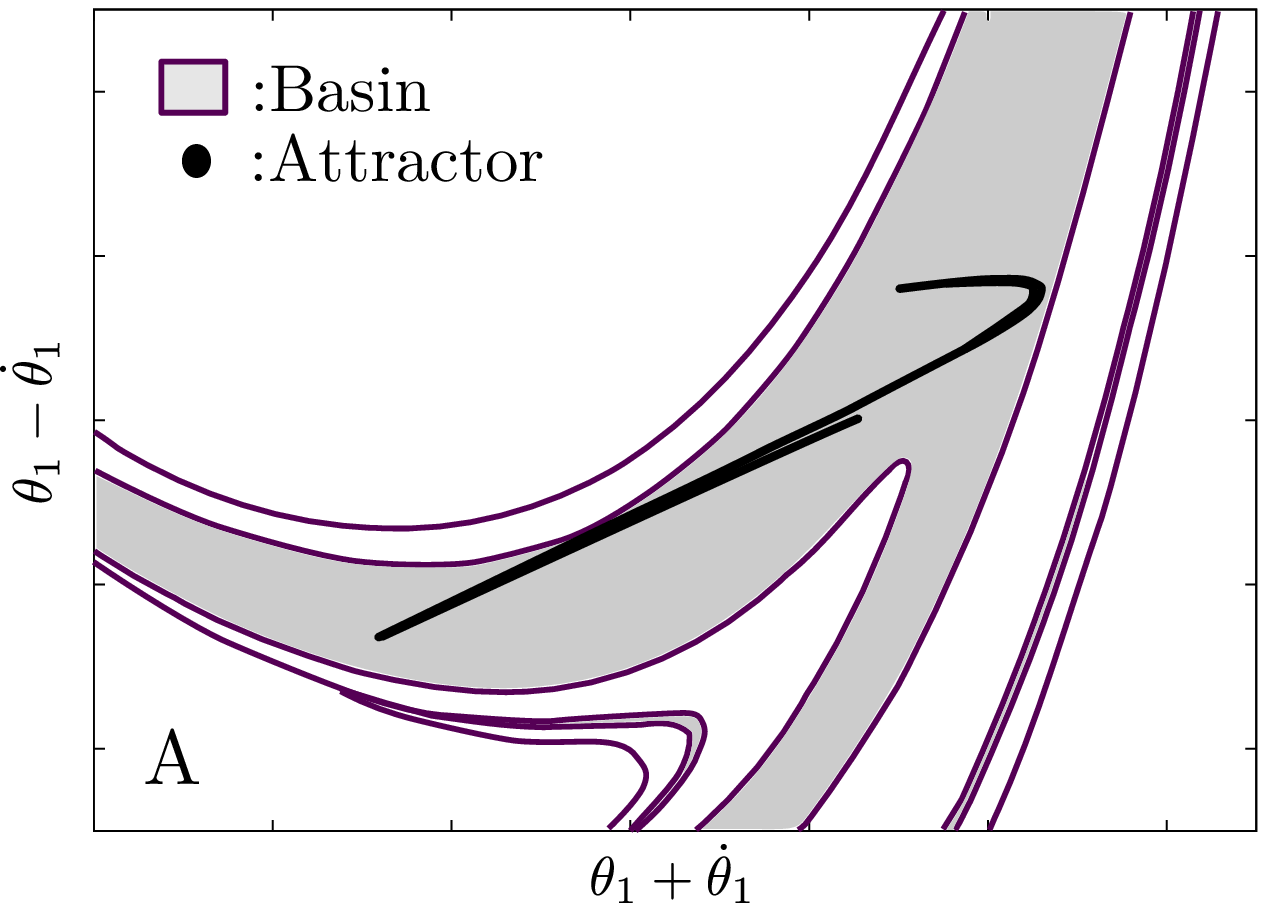} \ \ 
  \includegraphics[width=0.49\textwidth]{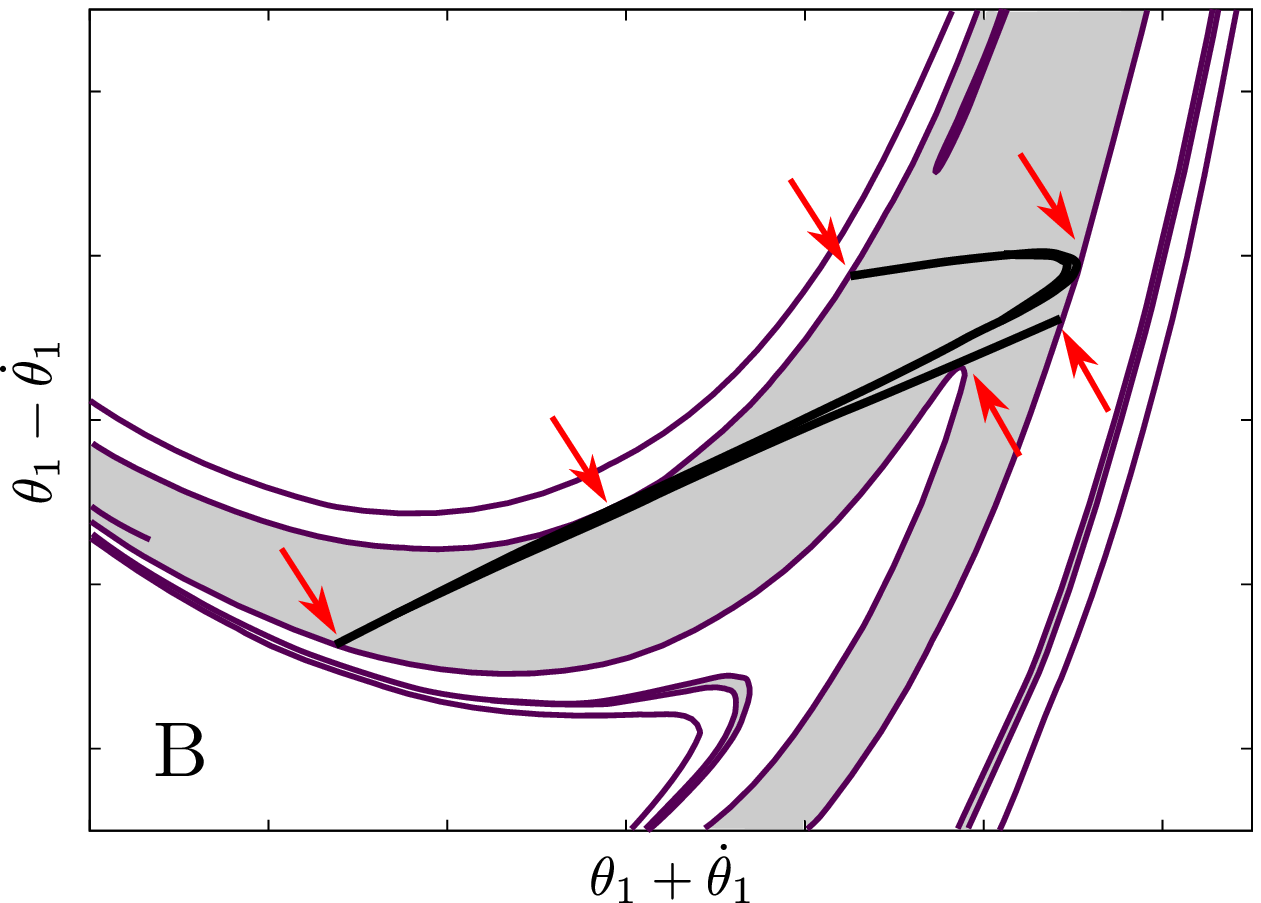}
  \caption{Geometric relationship between the chaotic attractor and the basin of attraction for $\gamma = 0.0187$ (A) and $0.01903$ (B)}\label{fig:contact_attr}
\end{figure}

In the above discussion, we clarified the mechanism that shapes the basin of attraction.
Although the parameter $\gamma$ induces a period-doubling cascade to chaos, the mechanism for constructing the basin of attraction is the same and does not depend on $\gamma$.
The domain, the inverse images of the domain, and the basin of attraction have only slight continuous changes, as shown in Fig.~\ref{fig:params}.
Figures~\ref{fig:contact_attr}A and B show the basin of attraction and the chaotic attractor for $\gamma=0.0187$ and $0.01903$, respectively, before the attractor disappears.
As $\gamma$ increases, the attractor expands to reach the boundaries of the basin of attraction.
When the attractor contacts the boundaries of the basin of attraction (red arrows in Fig.~\ref{fig:contact_attr}B), the attractor breaks down, and the model no longer continues walking.
In the theory of dynamical systems, this phenomenon is known as the boundary (attractor) crisis~\cite{Grebogi1983}.

\section{Conclusion}

In the present study, we clarified the construction mechanism for the basin of attraction for passive dynamic walking by focusing on the intrinsic hyperbolicity in the governing dynamics and using the viewpoint of the theory of dynamical systems.
We used the foot contact section $H$, the jump map $T$, the domain $D$ (the region on $T(H)$ where the model takes at least one step), the center-stable and center-unstable manifolds, the Poincar\'{e} map $S$, and the inverse maps.
Our results showed that the basin of attraction is determined by the relative positions of the center-unstable manifold and $T^{-1}(D)$, and the hyperbolicity near the saddle.
The equilibrium point and the manifolds are determined only by the continuous equations of motion, and they are independent of the foot-contact condition.
On the other hand, the positions of the domain and section are determined by the foot-contact condition and the jump map.
These inherent hybrid dynamics clarified the mechanism for constructing the basin of attraction, and concepts from the theory of dynamical systems, such as the center-stable and center-unstable manifolds, are very useful for the analysis of dynamic walking.

The thin, fractal-like basin of attraction of the simplest walking model is closely related to the one-dimensional instability of the upright equilibrium.
Because an inverted pendulum is governed by such a saddle-type instability, it plays important roles in the generation of various whole-body movements, such as body sway during quiet standing\cite{Suzuki2012, Asai2009, Funato2016}, as well as bipedal walking.
Although the present study focused on passive dynamic walking, our result is not specific to it, but is widely applicable to general bipedal walking, due to the intrinsic saddle property.

However, we note that the V-shaped basin of attraction, the slits, and the stripe patterns are formed by the relative positions of the center-unstable manifold and $T^{-1}(D)$, and the hyperbolicity near the saddle, as shown in Fig.~\ref{fig:domain_construction_1}.
Therefore, different bipedal walking models may have different shapes for the basin of attraction, depending on the relative positions of these regions.
However, due to the intrinsic saddle-type hyperbolicity, these elements have similar properties among bipedal walking models, so the discussion for these models may proceed in a similar way to those for our model.
Therefore, the present study may contribute not only to elucidating the stability mechanism in passive dynamic walking, but also to improving the understanding of the stability mechanism in human walking and to producing design principles for the control of walking support systems and biped robots.
In our future study, based on the geometrical characteristics clarified in this paper, we intend to improve the stability of bipedal walking by manipulating the relative positions of the center-unstable manifold and $T^{-1}(D)$ by designing a control system for a passive dynamic walking model.

\section*{Acknowledgements}

This paper is supported in part by Grant-in-Aid for Scientific Research (B) 15KT0015 from the Ministry of Education, Culture, Sports, Science, and Technology (MEXT) of Japan.
The main part of this research was done while the authors were supported by the JST CREST Project ``Alliance For Breakthrough Between Mathematics And Sciences''.

\bibliographystyle{ieeetr}
\bibliography{references}

\end{document}

%% file: basin_domain_1z.tex
\begingroup
  \makeatletter
  \providecommand\color[2][]{%
    \GenericError{(gnuplot) \space\space\space\@spaces}{%
      Package color not loaded in conjunction with
      terminal option `colourtext'%
    }{See the gnuplot documentation for explanation.%
    }{Either use 'blacktext' in gnuplot or load the package
      color.sty in LaTeX.}%
    \renewcommand\color[2][]{}%
  }%
  \providecommand\includegraphics[2][]{%
    \GenericError{(gnuplot) \space\space\space\@spaces}{%
      Package graphicx or graphics not loaded%
    }{See the gnuplot documentation for explanation.%
    }{The gnuplot epslatex terminal needs graphicx.sty or graphics.sty.}%
    \renewcommand\includegraphics[2][]{}%
  }%
  \providecommand\rotatebox[2]{#2}%
  \@ifundefined{ifGPcolor}{%
    \newif\ifGPcolor
    \GPcolortrue
  }{}%
  \@ifundefined{ifGPblacktext}{%
    \newif\ifGPblacktext
    \GPblacktexttrue
  }{}%
  \let\gplgaddtomacro\g@addto@macro
  \gdef\gplbacktext{}%
  \gdef\gplfronttext{}%
  \makeatother
  \ifGPblacktext
    \def\colorrgb#1{}%
    \def\colorgray#1{}%
  \else
    \ifGPcolor
      \def\colorrgb#1{\color[rgb]{#1}}%
      \def\colorgray#1{\color[gray]{#1}}%
      \expandafter\def\csname LTw\endcsname{\color{white}}%
      \expandafter\def\csname LTb\endcsname{\color{black}}%
      \expandafter\def\csname LTa\endcsname{\color{black}}%
      \expandafter\def\csname LT0\endcsname{\color[rgb]{1,0,0}}%
      \expandafter\def\csname LT1\endcsname{\color[rgb]{0,1,0}}%
      \expandafter\def\csname LT2\endcsname{\color[rgb]{0,0,1}}%
      \expandafter\def\csname LT3\endcsname{\color[rgb]{1,0,1}}%
      \expandafter\def\csname LT4\endcsname{\color[rgb]{0,1,1}}%
      \expandafter\def\csname LT5\endcsname{\color[rgb]{1,1,0}}%
      \expandafter\def\csname LT6\endcsname{\color[rgb]{0,0,0}}%
      \expandafter\def\csname LT7\endcsname{\color[rgb]{1,0.3,0}}%
      \expandafter\def\csname LT8\endcsname{\color[rgb]{0.5,0.5,0.5}}%
    \else
      \def\colorrgb#1{\color{black}}%
      \def\colorgray#1{\color[gray]{#1}}%
      \expandafter\def\csname LTw\endcsname{\color{white}}%
      \expandafter\def\csname LTb\endcsname{\color{black}}%
      \expandafter\def\csname LTa\endcsname{\color{black}}%
      \expandafter\def\csname LT0\endcsname{\color{black}}%
      \expandafter\def\csname LT1\endcsname{\color{black}}%
      \expandafter\def\csname LT2\endcsname{\color{black}}%
      \expandafter\def\csname LT3\endcsname{\color{black}}%
      \expandafter\def\csname LT4\endcsname{\color{black}}%
      \expandafter\def\csname LT5\endcsname{\color{black}}%
      \expandafter\def\csname LT6\endcsname{\color{black}}%
      \expandafter\def\csname LT7\endcsname{\color{black}}%
      \expandafter\def\csname LT8\endcsname{\color{black}}%
    \fi
  \fi
  \setlength{\unitlength}{0.0500bp}%
  \begin{picture}(7200.00,5040.00)%
    \gplgaddtomacro\gplbacktext{%
      \csname LTb\endcsname%
      \put(1078,704){\makebox(0,0)[r]{\strut{}-1}}%
      \put(1078,1156){\makebox(0,0)[r]{\strut{}-0.9}}%
      \put(1078,1609){\makebox(0,0)[r]{\strut{}-0.8}}%
      \put(1078,2061){\makebox(0,0)[r]{\strut{}-0.7}}%
      \put(1078,2513){\makebox(0,0)[r]{\strut{}-0.6}}%
      \put(1078,2966){\makebox(0,0)[r]{\strut{}-0.5}}%
      \put(1078,3418){\makebox(0,0)[r]{\strut{}-0.4}}%
      \put(1078,3870){\makebox(0,0)[r]{\strut{}-0.3}}%
      \put(1078,4323){\makebox(0,0)[r]{\strut{}-0.2}}%
      \put(1078,4775){\makebox(0,0)[r]{\strut{}-0.1}}%
      \put(1210,484){\makebox(0,0){\strut{} 0.1}}%
      \put(1839,484){\makebox(0,0){\strut{} 0.2}}%
      \put(2468,484){\makebox(0,0){\strut{} 0.3}}%
      \put(3096,484){\makebox(0,0){\strut{} 0.4}}%
      \put(3725,484){\makebox(0,0){\strut{} 0.5}}%
      \put(4354,484){\makebox(0,0){\strut{} 0.6}}%
      \put(4983,484){\makebox(0,0){\strut{} 0.7}}%
      \put(5611,484){\makebox(0,0){\strut{} 0.8}}%
      \put(6240,484){\makebox(0,0){\strut{} 0.9}}%
      \put(6869,484){\makebox(0,0){\strut{} 1}}%
      \put(308,2739){\rotatebox{-270}{\makebox(0,0){\strut{}\Large $\dot{\theta}_1$}}}%
      \put(4039,154){\makebox(0,0){\strut{}\Large $\theta_1$}}%
    }%
    \gplgaddtomacro\gplfronttext{%
    }%
    \gplbacktext
    \put(1240,700){\includegraphics{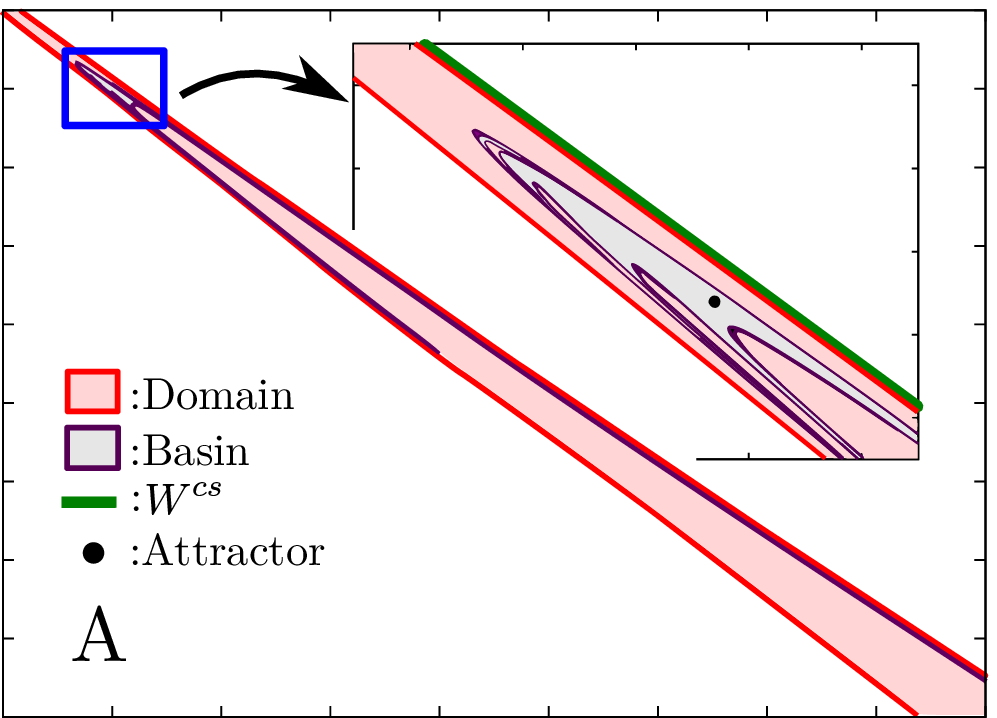}}%
    \gplfronttext
  \end{picture}%
\endgroup

%% file: basin_domain_3x.tex
\begingroup
  \makeatletter
  \providecommand\color[2][]{%
    \GenericError{(gnuplot) \space\space\space\@spaces}{%
      Package color not loaded in conjunction with
      terminal option `colourtext'%
    }{See the gnuplot documentation for explanation.%
    }{Either use 'blacktext' in gnuplot or load the package
      color.sty in LaTeX.}%
    \renewcommand\color[2][]{}%
  }%
  \providecommand\includegraphics[2][]{%
    \GenericError{(gnuplot) \space\space\space\@spaces}{%
      Package graphicx or graphics not loaded%
    }{See the gnuplot documentation for explanation.%
    }{The gnuplot epslatex terminal needs graphicx.sty or graphics.sty.}%
    \renewcommand\includegraphics[2][]{}%
  }%
  \providecommand\rotatebox[2]{#2}%
  \@ifundefined{ifGPcolor}{%
    \newif\ifGPcolor
    \GPcolortrue
  }{}%
  \@ifundefined{ifGPblacktext}{%
    \newif\ifGPblacktext
    \GPblacktexttrue
  }{}%
  \let\gplgaddtomacro\g@addto@macro
  \gdef\gplbacktext{}%
  \gdef\gplfronttext{}%
  \makeatother
  \ifGPblacktext
    \def\colorrgb#1{}%
    \def\colorgray#1{}%
  \else
    \ifGPcolor
      \def\colorrgb#1{\color[rgb]{#1}}%
      \def\colorgray#1{\color[gray]{#1}}%
      \expandafter\def\csname LTw\endcsname{\color{white}}%
      \expandafter\def\csname LTb\endcsname{\color{black}}%
      \expandafter\def\csname LTa\endcsname{\color{black}}%
      \expandafter\def\csname LT0\endcsname{\color[rgb]{1,0,0}}%
      \expandafter\def\csname LT1\endcsname{\color[rgb]{0,1,0}}%
      \expandafter\def\csname LT2\endcsname{\color[rgb]{0,0,1}}%
      \expandafter\def\csname LT3\endcsname{\color[rgb]{1,0,1}}%
      \expandafter\def\csname LT4\endcsname{\color[rgb]{0,1,1}}%
      \expandafter\def\csname LT5\endcsname{\color[rgb]{1,1,0}}%
      \expandafter\def\csname LT6\endcsname{\color[rgb]{0,0,0}}%
      \expandafter\def\csname LT7\endcsname{\color[rgb]{1,0.3,0}}%
      \expandafter\def\csname LT8\endcsname{\color[rgb]{0.5,0.5,0.5}}%
    \else
      \def\colorrgb#1{\color{black}}%
      \def\colorgray#1{\color[gray]{#1}}%
      \expandafter\def\csname LTw\endcsname{\color{white}}%
      \expandafter\def\csname LTb\endcsname{\color{black}}%
      \expandafter\def\csname LTa\endcsname{\color{black}}%
      \expandafter\def\csname LT0\endcsname{\color{black}}%
      \expandafter\def\csname LT1\endcsname{\color{black}}%
      \expandafter\def\csname LT2\endcsname{\color{black}}%
      \expandafter\def\csname LT3\endcsname{\color{black}}%
      \expandafter\def\csname LT4\endcsname{\color{black}}%
      \expandafter\def\csname LT5\endcsname{\color{black}}%
      \expandafter\def\csname LT6\endcsname{\color{black}}%
      \expandafter\def\csname LT7\endcsname{\color{black}}%
      \expandafter\def\csname LT8\endcsname{\color{black}}%
    \fi
  \fi
  \setlength{\unitlength}{0.0500bp}%
  \begin{picture}(7200.00,5040.00)%
    \gplgaddtomacro\gplbacktext{%
      \csname LTb\endcsname%
      \put(1078,704){\makebox(0,0)[r]{\strut{} 0.2}}%
      \put(1078,1213){\makebox(0,0)[r]{\strut{} 0.3}}%
      \put(1078,1722){\makebox(0,0)[r]{\strut{} 0.4}}%
      \put(1078,2231){\makebox(0,0)[r]{\strut{} 0.5}}%
      \put(1078,2740){\makebox(0,0)[r]{\strut{} 0.6}}%
      \put(1078,3248){\makebox(0,0)[r]{\strut{} 0.7}}%
      \put(1078,3757){\makebox(0,0)[r]{\strut{} 0.8}}%
      \put(1078,4266){\makebox(0,0)[r]{\strut{} 0.9}}%
      \put(1078,4775){\makebox(0,0)[r]{\strut{} 1}}%
      \put(1210,484){\makebox(0,0){\strut{}-0.04}}%
      \put(2153,484){\makebox(0,0){\strut{}-0.03}}%
      \put(3096,484){\makebox(0,0){\strut{}-0.02}}%
      \put(4040,484){\makebox(0,0){\strut{}-0.01}}%
      \put(4983,484){\makebox(0,0){\strut{} 0}}%
      \put(5926,484){\makebox(0,0){\strut{} 0.01}}%
      \put(6869,484){\makebox(0,0){\strut{} 0.02}}%
      \put(308,2739){\rotatebox{-270}{\makebox(0,0){\strut{}\Large $\theta_1 - \dot{\theta}_1$}}}%
      \put(4039,154){\makebox(0,0){\strut{}\Large $\theta_1 + \dot{\theta}_1$}}%
    }%
    \gplgaddtomacro\gplfronttext{%
    }%
    \gplbacktext
    \put(1240,700){\includegraphics{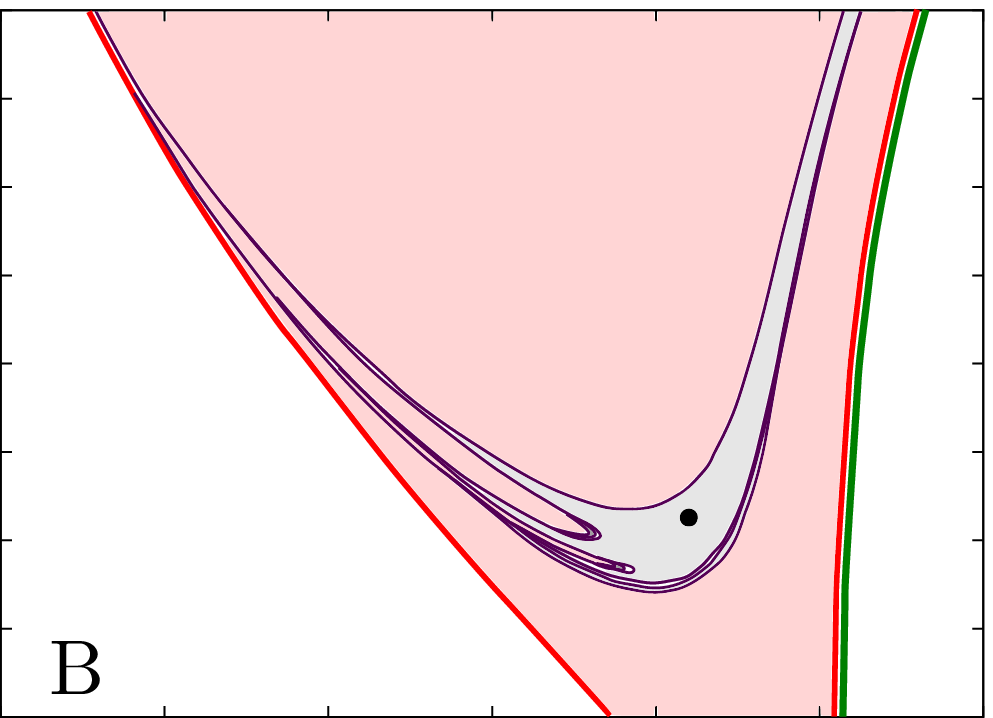}}%
    \gplfronttext
  \end{picture}%
\endgroup

%% file: domain123.tex
\begingroup
  \makeatletter
  \providecommand\color[2][]{%
    \GenericError{(gnuplot) \space\space\space\@spaces}{%
      Package color not loaded in conjunction with
      terminal option `colourtext'%
    }{See the gnuplot documentation for explanation.%
    }{Either use 'blacktext' in gnuplot or load the package
      color.sty in LaTeX.}%
    \renewcommand\color[2][]{}%
  }%
  \providecommand\includegraphics[2][]{%
    \GenericError{(gnuplot) \space\space\space\@spaces}{%
      Package graphicx or graphics not loaded%
    }{See the gnuplot documentation for explanation.%
    }{The gnuplot epslatex terminal needs graphicx.sty or graphics.sty.}%
    \renewcommand\includegraphics[2][]{}%
  }%
  \providecommand\rotatebox[2]{#2}%
  \@ifundefined{ifGPcolor}{%
    \newif\ifGPcolor
    \GPcolortrue
  }{}%
  \@ifundefined{ifGPblacktext}{%
    \newif\ifGPblacktext
    \GPblacktexttrue
  }{}%
  \let\gplgaddtomacro\g@addto@macro
  \gdef\gplbacktext{}%
  \gdef\gplfronttext{}%
  \makeatother
  \ifGPblacktext
    \def\colorrgb#1{}%
    \def\colorgray#1{}%
  \else
    \ifGPcolor
      \def\colorrgb#1{\color[rgb]{#1}}%
      \def\colorgray#1{\color[gray]{#1}}%
      \expandafter\def\csname LTw\endcsname{\color{white}}%
      \expandafter\def\csname LTb\endcsname{\color{black}}%
      \expandafter\def\csname LTa\endcsname{\color{black}}%
      \expandafter\def\csname LT0\endcsname{\color[rgb]{1,0,0}}%
      \expandafter\def\csname LT1\endcsname{\color[rgb]{0,1,0}}%
      \expandafter\def\csname LT2\endcsname{\color[rgb]{0,0,1}}%
      \expandafter\def\csname LT3\endcsname{\color[rgb]{1,0,1}}%
      \expandafter\def\csname LT4\endcsname{\color[rgb]{0,1,1}}%
      \expandafter\def\csname LT5\endcsname{\color[rgb]{1,1,0}}%
      \expandafter\def\csname LT6\endcsname{\color[rgb]{0,0,0}}%
      \expandafter\def\csname LT7\endcsname{\color[rgb]{1,0.3,0}}%
      \expandafter\def\csname LT8\endcsname{\color[rgb]{0.5,0.5,0.5}}%
    \else
      \def\colorrgb#1{\color{black}}%
      \def\colorgray#1{\color[gray]{#1}}%
      \expandafter\def\csname LTw\endcsname{\color{white}}%
      \expandafter\def\csname LTb\endcsname{\color{black}}%
      \expandafter\def\csname LTa\endcsname{\color{black}}%
      \expandafter\def\csname LT0\endcsname{\color{black}}%
      \expandafter\def\csname LT1\endcsname{\color{black}}%
      \expandafter\def\csname LT2\endcsname{\color{black}}%
      \expandafter\def\csname LT3\endcsname{\color{black}}%
      \expandafter\def\csname LT4\endcsname{\color{black}}%
      \expandafter\def\csname LT5\endcsname{\color{black}}%
      \expandafter\def\csname LT6\endcsname{\color{black}}%
      \expandafter\def\csname LT7\endcsname{\color{black}}%
      \expandafter\def\csname LT8\endcsname{\color{black}}%
    \fi
  \fi
  \setlength{\unitlength}{0.0500bp}%
  \begin{picture}(7200.00,5040.00)%
    \gplgaddtomacro\gplbacktext{%
      \csname LTb\endcsname%
      \put(1078,704){\makebox(0,0)[r]{\strut{} 0.2}}%
      \put(1078,1213){\makebox(0,0)[r]{\strut{} 0.3}}%
      \put(1078,1722){\makebox(0,0)[r]{\strut{} 0.4}}%
      \put(1078,2231){\makebox(0,0)[r]{\strut{} 0.5}}%
      \put(1078,2740){\makebox(0,0)[r]{\strut{} 0.6}}%
      \put(1078,3248){\makebox(0,0)[r]{\strut{} 0.7}}%
      \put(1078,3757){\makebox(0,0)[r]{\strut{} 0.8}}%
      \put(1078,4266){\makebox(0,0)[r]{\strut{} 0.9}}%
      \put(1078,4775){\makebox(0,0)[r]{\strut{} 1}}%
      \put(1210,484){\makebox(0,0){\strut{}-0.04}}%
      \put(2153,484){\makebox(0,0){\strut{}-0.03}}%
      \put(3096,484){\makebox(0,0){\strut{}-0.02}}%
      \put(4040,484){\makebox(0,0){\strut{}-0.01}}%
      \put(4983,484){\makebox(0,0){\strut{} 0}}%
      \put(5926,484){\makebox(0,0){\strut{} 0.01}}%
      \put(6869,484){\makebox(0,0){\strut{} 0.02}}%
      \put(308,2739){\rotatebox{-270}{\makebox(0,0){\strut{}{\Large $\theta_1 - \dot{\theta}_1$}}}}%
      \put(4039,154){\makebox(0,0){\strut{}{\Large $\theta_1 + \dot{\theta}_1$}}}%
    }%
    \gplgaddtomacro\gplfronttext{%
    }%
    \gplbacktext
    \put(1240,700){\includegraphics{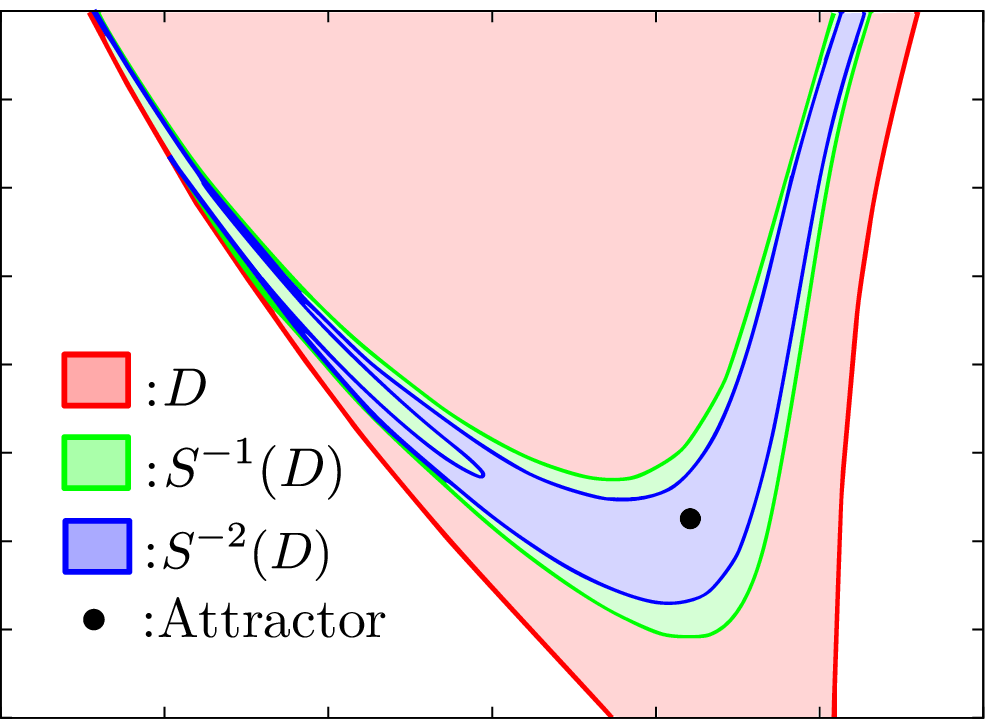}}%
    \gplfronttext
  \end{picture}%
\endgroup

%% file: domain8_8L.tex
\begingroup
  \makeatletter
  \providecommand\color[2][]{%
    \GenericError{(gnuplot) \space\space\space\@spaces}{%
      Package color not loaded in conjunction with
      terminal option `colourtext'%
    }{See the gnuplot documentation for explanation.%
    }{Either use 'blacktext' in gnuplot or load the package
      color.sty in LaTeX.}%
    \renewcommand\color[2][]{}%
  }%
  \providecommand\includegraphics[2][]{%
    \GenericError{(gnuplot) \space\space\space\@spaces}{%
      Package graphicx or graphics not loaded%
    }{See the gnuplot documentation for explanation.%
    }{The gnuplot epslatex terminal needs graphicx.sty or graphics.sty.}%
    \renewcommand\includegraphics[2][]{}%
  }%
  \providecommand\rotatebox[2]{#2}%
  \@ifundefined{ifGPcolor}{%
    \newif\ifGPcolor
    \GPcolortrue
  }{}%
  \@ifundefined{ifGPblacktext}{%
    \newif\ifGPblacktext
    \GPblacktexttrue
  }{}%
  \let\gplgaddtomacro\g@addto@macro
  \gdef\gplbacktext{}%
  \gdef\gplfronttext{}%
  \makeatother
  \ifGPblacktext
    \def\colorrgb#1{}%
    \def\colorgray#1{}%
  \else
    \ifGPcolor
      \def\colorrgb#1{\color[rgb]{#1}}%
      \def\colorgray#1{\color[gray]{#1}}%
      \expandafter\def\csname LTw\endcsname{\color{white}}%
      \expandafter\def\csname LTb\endcsname{\color{black}}%
      \expandafter\def\csname LTa\endcsname{\color{black}}%
      \expandafter\def\csname LT0\endcsname{\color[rgb]{1,0,0}}%
      \expandafter\def\csname LT1\endcsname{\color[rgb]{0,1,0}}%
      \expandafter\def\csname LT2\endcsname{\color[rgb]{0,0,1}}%
      \expandafter\def\csname LT3\endcsname{\color[rgb]{1,0,1}}%
      \expandafter\def\csname LT4\endcsname{\color[rgb]{0,1,1}}%
      \expandafter\def\csname LT5\endcsname{\color[rgb]{1,1,0}}%
      \expandafter\def\csname LT6\endcsname{\color[rgb]{0,0,0}}%
      \expandafter\def\csname LT7\endcsname{\color[rgb]{1,0.3,0}}%
      \expandafter\def\csname LT8\endcsname{\color[rgb]{0.5,0.5,0.5}}%
    \else
      \def\colorrgb#1{\color{black}}%
      \def\colorgray#1{\color[gray]{#1}}%
      \expandafter\def\csname LTw\endcsname{\color{white}}%
      \expandafter\def\csname LTb\endcsname{\color{black}}%
      \expandafter\def\csname LTa\endcsname{\color{black}}%
      \expandafter\def\csname LT0\endcsname{\color{black}}%
      \expandafter\def\csname LT1\endcsname{\color{black}}%
      \expandafter\def\csname LT2\endcsname{\color{black}}%
      \expandafter\def\csname LT3\endcsname{\color{black}}%
      \expandafter\def\csname LT4\endcsname{\color{black}}%
      \expandafter\def\csname LT5\endcsname{\color{black}}%
      \expandafter\def\csname LT6\endcsname{\color{black}}%
      \expandafter\def\csname LT7\endcsname{\color{black}}%
      \expandafter\def\csname LT8\endcsname{\color{black}}%
    \fi
  \fi
  \setlength{\unitlength}{0.0500bp}%
  \begin{picture}(7200.00,5040.00)%
    \gplgaddtomacro\gplbacktext{%
      \csname LTb\endcsname%
      \put(1078,704){\makebox(0,0)[r]{\strut{} 0.3}}%
      \put(1078,1286){\makebox(0,0)[r]{\strut{} 0.4}}%
      \put(1078,1867){\makebox(0,0)[r]{\strut{} 0.5}}%
      \put(1078,2449){\makebox(0,0)[r]{\strut{} 0.6}}%
      \put(1078,3030){\makebox(0,0)[r]{\strut{} 0.7}}%
      \put(1078,3612){\makebox(0,0)[r]{\strut{} 0.8}}%
      \put(1078,4193){\makebox(0,0)[r]{\strut{} 0.9}}%
      \put(1078,4775){\makebox(0,0)[r]{\strut{} 1}}%
      \put(1210,484){\makebox(0,0){\strut{}-0.035}}%
      \put(1776,484){\makebox(0,0){\strut{}-0.03}}%
      \put(2342,484){\makebox(0,0){\strut{}-0.025}}%
      \put(2908,484){\makebox(0,0){\strut{}-0.02}}%
      \put(3474,484){\makebox(0,0){\strut{}-0.015}}%
      \put(4040,484){\makebox(0,0){\strut{}-0.01}}%
      \put(4605,484){\makebox(0,0){\strut{}-0.005}}%
      \put(5171,484){\makebox(0,0){\strut{} 0}}%
      \put(5737,484){\makebox(0,0){\strut{} 0.005}}%
      \put(6303,484){\makebox(0,0){\strut{} 0.01}}%
      \put(6869,484){\makebox(0,0){\strut{} 0.015}}%
      \put(308,2739){\rotatebox{-270}{\makebox(0,0){\strut{}\Large $\theta_1 - \dot{\theta}_1$}}}%
      \put(4039,154){\makebox(0,0){\strut{}\Large $\theta_1 + \dot{\theta}_1$}}%
    }%
    \gplgaddtomacro\gplfronttext{%
    }%
    \gplbacktext
    \put(1240,700){\includegraphics{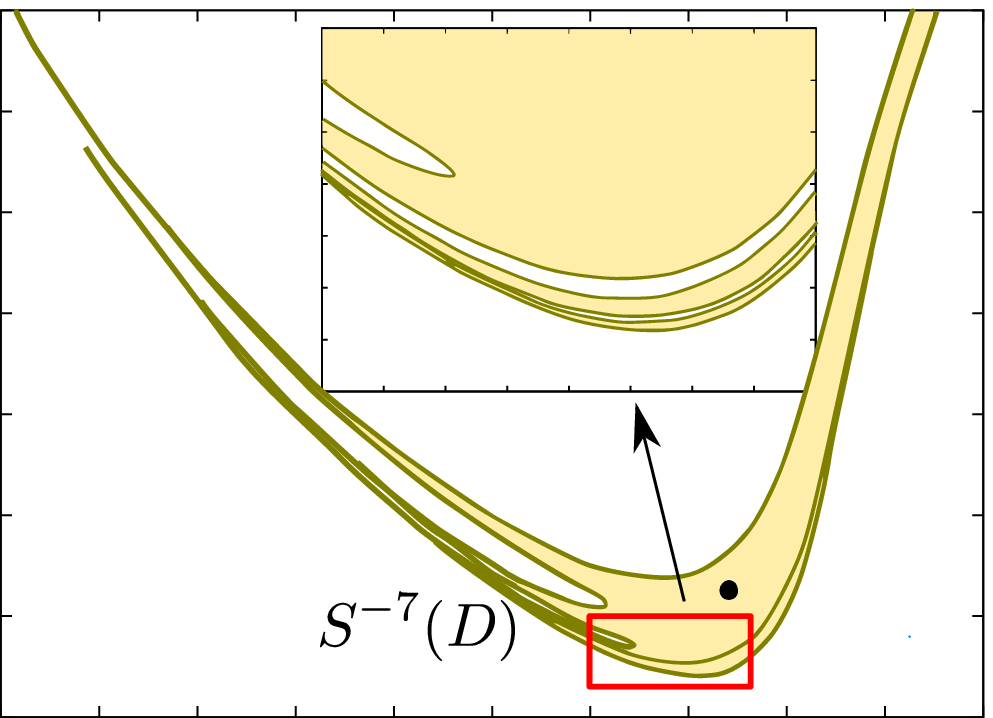}}%
    \gplfronttext
  \end{picture}%
\endgroup

%% file: 0_0160_basin.tex
\begingroup
  \makeatletter
  \providecommand\color[2][]{%
    \GenericError{(gnuplot) \space\space\space\@spaces}{%
      Package color not loaded in conjunction with
      terminal option `colourtext'%
    }{See the gnuplot documentation for explanation.%
    }{Either use 'blacktext' in gnuplot or load the package
      color.sty in LaTeX.}%
    \renewcommand\color[2][]{}%
  }%
  \providecommand\includegraphics[2][]{%
    \GenericError{(gnuplot) \space\space\space\@spaces}{%
      Package graphicx or graphics not loaded%
    }{See the gnuplot documentation for explanation.%
    }{The gnuplot epslatex terminal needs graphicx.sty or graphics.sty.}%
    \renewcommand\includegraphics[2][]{}%
  }%
  \providecommand\rotatebox[2]{#2}%
  \@ifundefined{ifGPcolor}{%
    \newif\ifGPcolor
    \GPcolortrue
  }{}%
  \@ifundefined{ifGPblacktext}{%
    \newif\ifGPblacktext
    \GPblacktexttrue
  }{}%
  \let\gplgaddtomacro\g@addto@macro
  \gdef\gplbacktext{}%
  \gdef\gplfronttext{}%
  \makeatother
  \ifGPblacktext
    \def\colorrgb#1{}%
    \def\colorgray#1{}%
  \else
    \ifGPcolor
      \def\colorrgb#1{\color[rgb]{#1}}%
      \def\colorgray#1{\color[gray]{#1}}%
      \expandafter\def\csname LTw\endcsname{\color{white}}%
      \expandafter\def\csname LTb\endcsname{\color{black}}%
      \expandafter\def\csname LTa\endcsname{\color{black}}%
      \expandafter\def\csname LT0\endcsname{\color[rgb]{1,0,0}}%
      \expandafter\def\csname LT1\endcsname{\color[rgb]{0,1,0}}%
      \expandafter\def\csname LT2\endcsname{\color[rgb]{0,0,1}}%
      \expandafter\def\csname LT3\endcsname{\color[rgb]{1,0,1}}%
      \expandafter\def\csname LT4\endcsname{\color[rgb]{0,1,1}}%
      \expandafter\def\csname LT5\endcsname{\color[rgb]{1,1,0}}%
      \expandafter\def\csname LT6\endcsname{\color[rgb]{0,0,0}}%
      \expandafter\def\csname LT7\endcsname{\color[rgb]{1,0.3,0}}%
      \expandafter\def\csname LT8\endcsname{\color[rgb]{0.5,0.5,0.5}}%
    \else
      \def\colorrgb#1{\color{black}}%
      \def\colorgray#1{\color[gray]{#1}}%
      \expandafter\def\csname LTw\endcsname{\color{white}}%
      \expandafter\def\csname LTb\endcsname{\color{black}}%
      \expandafter\def\csname LTa\endcsname{\color{black}}%
      \expandafter\def\csname LT0\endcsname{\color{black}}%
      \expandafter\def\csname LT1\endcsname{\color{black}}%
      \expandafter\def\csname LT2\endcsname{\color{black}}%
      \expandafter\def\csname LT3\endcsname{\color{black}}%
      \expandafter\def\csname LT4\endcsname{\color{black}}%
      \expandafter\def\csname LT5\endcsname{\color{black}}%
      \expandafter\def\csname LT6\endcsname{\color{black}}%
      \expandafter\def\csname LT7\endcsname{\color{black}}%
      \expandafter\def\csname LT8\endcsname{\color{black}}%
    \fi
  \fi
  \setlength{\unitlength}{0.0500bp}%
  \begin{picture}(7200.00,5040.00)%
    \gplgaddtomacro\gplbacktext{%
      \csname LTb\endcsname%
      \put(1078,704){\makebox(0,0)[r]{\strut{} 0.2}}%
      \put(1078,1213){\makebox(0,0)[r]{\strut{} 0.3}}%
      \put(1078,1722){\makebox(0,0)[r]{\strut{} 0.4}}%
      \put(1078,2231){\makebox(0,0)[r]{\strut{} 0.5}}%
      \put(1078,2740){\makebox(0,0)[r]{\strut{} 0.6}}%
      \put(1078,3248){\makebox(0,0)[r]{\strut{} 0.7}}%
      \put(1078,3757){\makebox(0,0)[r]{\strut{} 0.8}}%
      \put(1078,4266){\makebox(0,0)[r]{\strut{} 0.9}}%
      \put(1078,4775){\makebox(0,0)[r]{\strut{} 1}}%
      \put(1210,484){\makebox(0,0){\strut{}-0.04}}%
      \put(2153,484){\makebox(0,0){\strut{}-0.03}}%
      \put(3096,484){\makebox(0,0){\strut{}-0.02}}%
      \put(4040,484){\makebox(0,0){\strut{}-0.01}}%
      \put(4983,484){\makebox(0,0){\strut{} 0}}%
      \put(5926,484){\makebox(0,0){\strut{} 0.01}}%
      \put(6869,484){\makebox(0,0){\strut{} 0.02}}%
      \put(308,2739){\rotatebox{-270}{\makebox(0,0){\strut{}\Large $\theta_1 - \dot{\theta}_1$}}}%
      \put(4039,154){\makebox(0,0){\strut{}\Large $\theta_1 + \dot{\theta}_1$}}%
    }%
    \gplgaddtomacro\gplfronttext{%
    }%
    \gplbacktext
    \put(1240,700){\includegraphics{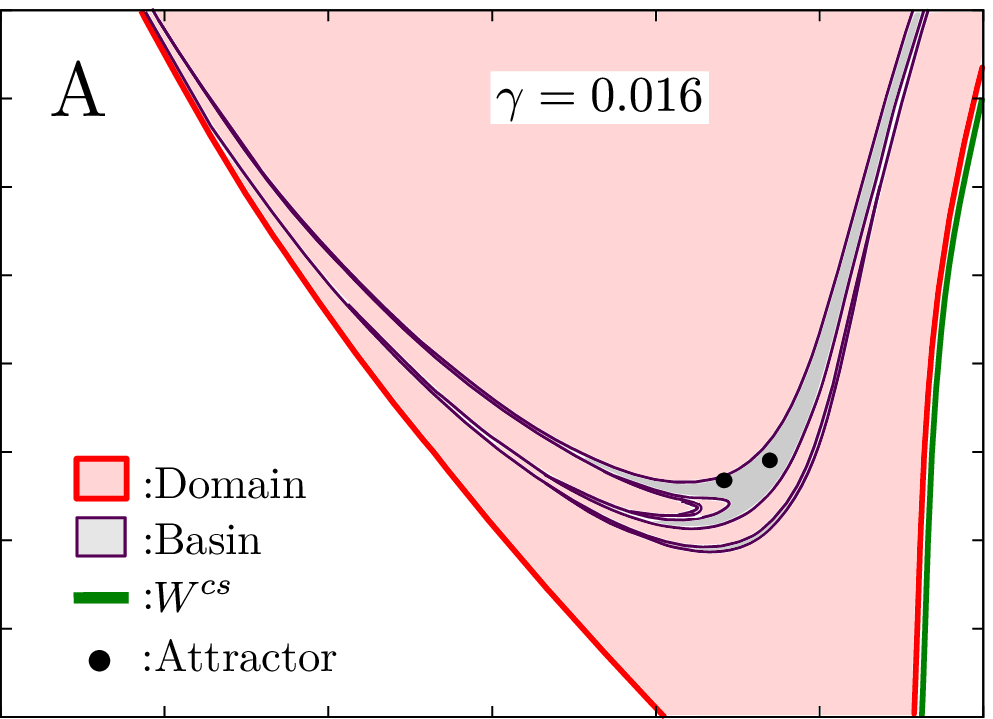}}%
    \gplfronttext
  \end{picture}%
\endgroup

%% file: 0_0160.tex
\begingroup
  \makeatletter
  \providecommand\color[2][]{%
    \GenericError{(gnuplot) \space\space\space\@spaces}{%
      Package color not loaded in conjunction with
      terminal option `colourtext'%
    }{See the gnuplot documentation for explanation.%
    }{Either use 'blacktext' in gnuplot or load the package
      color.sty in LaTeX.}%
    \renewcommand\color[2][]{}%
  }%
  \providecommand\includegraphics[2][]{%
    \GenericError{(gnuplot) \space\space\space\@spaces}{%
      Package graphicx or graphics not loaded%
    }{See the gnuplot documentation for explanation.%
    }{The gnuplot epslatex terminal needs graphicx.sty or graphics.sty.}%
    \renewcommand\includegraphics[2][]{}%
  }%
  \providecommand\rotatebox[2]{#2}%
  \@ifundefined{ifGPcolor}{%
    \newif\ifGPcolor
    \GPcolortrue
  }{}%
  \@ifundefined{ifGPblacktext}{%
    \newif\ifGPblacktext
    \GPblacktexttrue
  }{}%
  \let\gplgaddtomacro\g@addto@macro
  \gdef\gplbacktext{}%
  \gdef\gplfronttext{}%
  \makeatother
  \ifGPblacktext
    \def\colorrgb#1{}%
    \def\colorgray#1{}%
  \else
    \ifGPcolor
      \def\colorrgb#1{\color[rgb]{#1}}%
      \def\colorgray#1{\color[gray]{#1}}%
      \expandafter\def\csname LTw\endcsname{\color{white}}%
      \expandafter\def\csname LTb\endcsname{\color{black}}%
      \expandafter\def\csname LTa\endcsname{\color{black}}%
      \expandafter\def\csname LT0\endcsname{\color[rgb]{1,0,0}}%
      \expandafter\def\csname LT1\endcsname{\color[rgb]{0,1,0}}%
      \expandafter\def\csname LT2\endcsname{\color[rgb]{0,0,1}}%
      \expandafter\def\csname LT3\endcsname{\color[rgb]{1,0,1}}%
      \expandafter\def\csname LT4\endcsname{\color[rgb]{0,1,1}}%
      \expandafter\def\csname LT5\endcsname{\color[rgb]{1,1,0}}%
      \expandafter\def\csname LT6\endcsname{\color[rgb]{0,0,0}}%
      \expandafter\def\csname LT7\endcsname{\color[rgb]{1,0.3,0}}%
      \expandafter\def\csname LT8\endcsname{\color[rgb]{0.5,0.5,0.5}}%
    \else
      \def\colorrgb#1{\color{black}}%
      \def\colorgray#1{\color[gray]{#1}}%
      \expandafter\def\csname LTw\endcsname{\color{white}}%
      \expandafter\def\csname LTb\endcsname{\color{black}}%
      \expandafter\def\csname LTa\endcsname{\color{black}}%
      \expandafter\def\csname LT0\endcsname{\color{black}}%
      \expandafter\def\csname LT1\endcsname{\color{black}}%
      \expandafter\def\csname LT2\endcsname{\color{black}}%
      \expandafter\def\csname LT3\endcsname{\color{black}}%
      \expandafter\def\csname LT4\endcsname{\color{black}}%
      \expandafter\def\csname LT5\endcsname{\color{black}}%
      \expandafter\def\csname LT6\endcsname{\color{black}}%
      \expandafter\def\csname LT7\endcsname{\color{black}}%
      \expandafter\def\csname LT8\endcsname{\color{black}}%
    \fi
  \fi
  \setlength{\unitlength}{0.0500bp}%
  \begin{picture}(7200.00,5040.00)%
    \gplgaddtomacro\gplbacktext{%
      \csname LTb\endcsname%
      \put(1078,704){\makebox(0,0)[r]{\strut{} 0.2}}%
      \put(1078,1213){\makebox(0,0)[r]{\strut{} 0.3}}%
      \put(1078,1722){\makebox(0,0)[r]{\strut{} 0.4}}%
      \put(1078,2231){\makebox(0,0)[r]{\strut{} 0.5}}%
      \put(1078,2740){\makebox(0,0)[r]{\strut{} 0.6}}%
      \put(1078,3248){\makebox(0,0)[r]{\strut{} 0.7}}%
      \put(1078,3757){\makebox(0,0)[r]{\strut{} 0.8}}%
      \put(1078,4266){\makebox(0,0)[r]{\strut{} 0.9}}%
      \put(1078,4775){\makebox(0,0)[r]{\strut{} 1}}%
      \put(1210,484){\makebox(0,0){\strut{}-0.04}}%
      \put(2153,484){\makebox(0,0){\strut{}-0.03}}%
      \put(3096,484){\makebox(0,0){\strut{}-0.02}}%
      \put(4040,484){\makebox(0,0){\strut{}-0.01}}%
      \put(4983,484){\makebox(0,0){\strut{} 0}}%
      \put(5926,484){\makebox(0,0){\strut{} 0.01}}%
      \put(6869,484){\makebox(0,0){\strut{} 0.02}}%
      \put(308,2739){\rotatebox{-270}{\makebox(0,0){\strut{} \Large $\theta_1 - \dot{\theta}_1$}}}%
      \put(4039,154){\makebox(0,0){\strut{} \Large $\theta_1 + \dot{\theta}_1$}}%
    }%
    \gplgaddtomacro\gplfronttext{%
    }%
    \gplbacktext
    \put(1240,700){\includegraphics{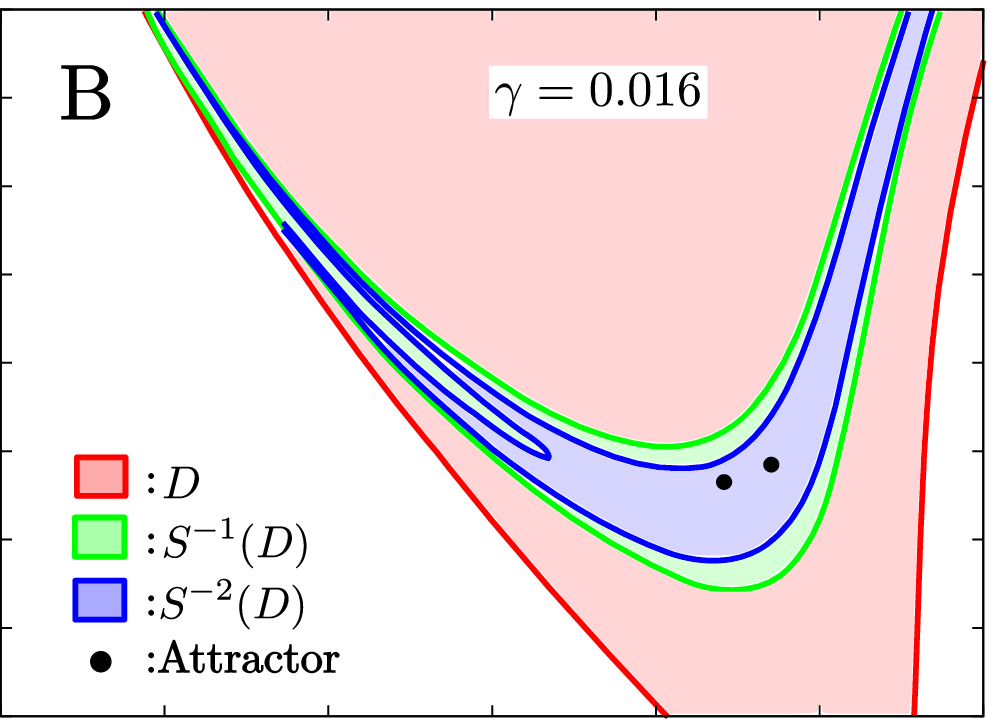}}%
    \gplfronttext
  \end{picture}%
\endgroup

%% file: 0_01780_basin.tex
\begingroup
  \makeatletter
  \providecommand\color[2][]{%
    \GenericError{(gnuplot) \space\space\space\@spaces}{%
      Package color not loaded in conjunction with
      terminal option `colourtext'%
    }{See the gnuplot documentation for explanation.%
    }{Either use 'blacktext' in gnuplot or load the package
      color.sty in LaTeX.}%
    \renewcommand\color[2][]{}%
  }%
  \providecommand\includegraphics[2][]{%
    \GenericError{(gnuplot) \space\space\space\@spaces}{%
      Package graphicx or graphics not loaded%
    }{See the gnuplot documentation for explanation.%
    }{The gnuplot epslatex terminal needs graphicx.sty or graphics.sty.}%
    \renewcommand\includegraphics[2][]{}%
  }%
  \providecommand\rotatebox[2]{#2}%
  \@ifundefined{ifGPcolor}{%
    \newif\ifGPcolor
    \GPcolortrue
  }{}%
  \@ifundefined{ifGPblacktext}{%
    \newif\ifGPblacktext
    \GPblacktexttrue
  }{}%
  \let\gplgaddtomacro\g@addto@macro
  \gdef\gplbacktext{}%
  \gdef\gplfronttext{}%
  \makeatother
  \ifGPblacktext
    \def\colorrgb#1{}%
    \def\colorgray#1{}%
  \else
    \ifGPcolor
      \def\colorrgb#1{\color[rgb]{#1}}%
      \def\colorgray#1{\color[gray]{#1}}%
      \expandafter\def\csname LTw\endcsname{\color{white}}%
      \expandafter\def\csname LTb\endcsname{\color{black}}%
      \expandafter\def\csname LTa\endcsname{\color{black}}%
      \expandafter\def\csname LT0\endcsname{\color[rgb]{1,0,0}}%
      \expandafter\def\csname LT1\endcsname{\color[rgb]{0,1,0}}%
      \expandafter\def\csname LT2\endcsname{\color[rgb]{0,0,1}}%
      \expandafter\def\csname LT3\endcsname{\color[rgb]{1,0,1}}%
      \expandafter\def\csname LT4\endcsname{\color[rgb]{0,1,1}}%
      \expandafter\def\csname LT5\endcsname{\color[rgb]{1,1,0}}%
      \expandafter\def\csname LT6\endcsname{\color[rgb]{0,0,0}}%
      \expandafter\def\csname LT7\endcsname{\color[rgb]{1,0.3,0}}%
      \expandafter\def\csname LT8\endcsname{\color[rgb]{0.5,0.5,0.5}}%
    \else
      \def\colorrgb#1{\color{black}}%
      \def\colorgray#1{\color[gray]{#1}}%
      \expandafter\def\csname LTw\endcsname{\color{white}}%
      \expandafter\def\csname LTb\endcsname{\color{black}}%
      \expandafter\def\csname LTa\endcsname{\color{black}}%
      \expandafter\def\csname LT0\endcsname{\color{black}}%
      \expandafter\def\csname LT1\endcsname{\color{black}}%
      \expandafter\def\csname LT2\endcsname{\color{black}}%
      \expandafter\def\csname LT3\endcsname{\color{black}}%
      \expandafter\def\csname LT4\endcsname{\color{black}}%
      \expandafter\def\csname LT5\endcsname{\color{black}}%
      \expandafter\def\csname LT6\endcsname{\color{black}}%
      \expandafter\def\csname LT7\endcsname{\color{black}}%
      \expandafter\def\csname LT8\endcsname{\color{black}}%
    \fi
  \fi
  \setlength{\unitlength}{0.0500bp}%
  \begin{picture}(7200.00,5040.00)%
    \gplgaddtomacro\gplbacktext{%
      \csname LTb\endcsname%
      \put(1078,704){\makebox(0,0)[r]{\strut{} 0.2}}%
      \put(1078,1213){\makebox(0,0)[r]{\strut{} 0.3}}%
      \put(1078,1722){\makebox(0,0)[r]{\strut{} 0.4}}%
      \put(1078,2231){\makebox(0,0)[r]{\strut{} 0.5}}%
      \put(1078,2740){\makebox(0,0)[r]{\strut{} 0.6}}%
      \put(1078,3248){\makebox(0,0)[r]{\strut{} 0.7}}%
      \put(1078,3757){\makebox(0,0)[r]{\strut{} 0.8}}%
      \put(1078,4266){\makebox(0,0)[r]{\strut{} 0.9}}%
      \put(1078,4775){\makebox(0,0)[r]{\strut{} 1}}%
      \put(1210,484){\makebox(0,0){\strut{}-0.04}}%
      \put(2153,484){\makebox(0,0){\strut{}-0.03}}%
      \put(3096,484){\makebox(0,0){\strut{}-0.02}}%
      \put(4040,484){\makebox(0,0){\strut{}-0.01}}%
      \put(4983,484){\makebox(0,0){\strut{} 0}}%
      \put(5926,484){\makebox(0,0){\strut{} 0.01}}%
      \put(6869,484){\makebox(0,0){\strut{} 0.02}}%
      \put(308,2739){\rotatebox{-270}{\makebox(0,0){\strut{}\Large $\theta_1 - \dot{\theta}_1$}}}%
      \put(4039,154){\makebox(0,0){\strut{}\Large $\theta_1 + \dot{\theta}_1$}}%
    }%
    \gplgaddtomacro\gplfronttext{%
    }%
    \gplbacktext
    \put(1240,700){\includegraphics{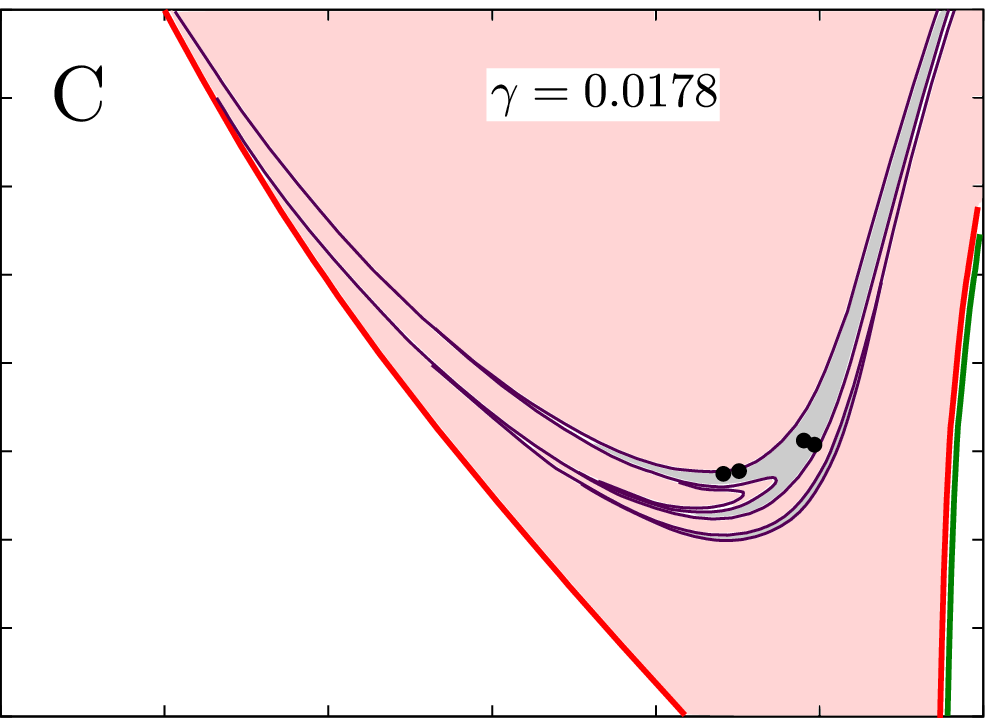}}%
    \gplfronttext
  \end{picture}%
\endgroup

%% file: 0_01780.tex
\begingroup
  \makeatletter
  \providecommand\color[2][]{%
    \GenericError{(gnuplot) \space\space\space\@spaces}{%
      Package color not loaded in conjunction with
      terminal option `colourtext'%
    }{See the gnuplot documentation for explanation.%
    }{Either use 'blacktext' in gnuplot or load the package
      color.sty in LaTeX.}%
    \renewcommand\color[2][]{}%
  }%
  \providecommand\includegraphics[2][]{%
    \GenericError{(gnuplot) \space\space\space\@spaces}{%
      Package graphicx or graphics not loaded%
    }{See the gnuplot documentation for explanation.%
    }{The gnuplot epslatex terminal needs graphicx.sty or graphics.sty.}%
    \renewcommand\includegraphics[2][]{}%
  }%
  \providecommand\rotatebox[2]{#2}%
  \@ifundefined{ifGPcolor}{%
    \newif\ifGPcolor
    \GPcolortrue
  }{}%
  \@ifundefined{ifGPblacktext}{%
    \newif\ifGPblacktext
    \GPblacktexttrue
  }{}%
  \let\gplgaddtomacro\g@addto@macro
  \gdef\gplbacktext{}%
  \gdef\gplfronttext{}%
  \makeatother
  \ifGPblacktext
    \def\colorrgb#1{}%
    \def\colorgray#1{}%
  \else
    \ifGPcolor
      \def\colorrgb#1{\color[rgb]{#1}}%
      \def\colorgray#1{\color[gray]{#1}}%
      \expandafter\def\csname LTw\endcsname{\color{white}}%
      \expandafter\def\csname LTb\endcsname{\color{black}}%
      \expandafter\def\csname LTa\endcsname{\color{black}}%
      \expandafter\def\csname LT0\endcsname{\color[rgb]{1,0,0}}%
      \expandafter\def\csname LT1\endcsname{\color[rgb]{0,1,0}}%
      \expandafter\def\csname LT2\endcsname{\color[rgb]{0,0,1}}%
      \expandafter\def\csname LT3\endcsname{\color[rgb]{1,0,1}}%
      \expandafter\def\csname LT4\endcsname{\color[rgb]{0,1,1}}%
      \expandafter\def\csname LT5\endcsname{\color[rgb]{1,1,0}}%
      \expandafter\def\csname LT6\endcsname{\color[rgb]{0,0,0}}%
      \expandafter\def\csname LT7\endcsname{\color[rgb]{1,0.3,0}}%
      \expandafter\def\csname LT8\endcsname{\color[rgb]{0.5,0.5,0.5}}%
    \else
      \def\colorrgb#1{\color{black}}%
      \def\colorgray#1{\color[gray]{#1}}%
      \expandafter\def\csname LTw\endcsname{\color{white}}%
      \expandafter\def\csname LTb\endcsname{\color{black}}%
      \expandafter\def\csname LTa\endcsname{\color{black}}%
      \expandafter\def\csname LT0\endcsname{\color{black}}%
      \expandafter\def\csname LT1\endcsname{\color{black}}%
      \expandafter\def\csname LT2\endcsname{\color{black}}%
      \expandafter\def\csname LT3\endcsname{\color{black}}%
      \expandafter\def\csname LT4\endcsname{\color{black}}%
      \expandafter\def\csname LT5\endcsname{\color{black}}%
      \expandafter\def\csname LT6\endcsname{\color{black}}%
      \expandafter\def\csname LT7\endcsname{\color{black}}%
      \expandafter\def\csname LT8\endcsname{\color{black}}%
    \fi
  \fi
  \setlength{\unitlength}{0.0500bp}%
  \begin{picture}(7200.00,5040.00)%
    \gplgaddtomacro\gplbacktext{%
      \csname LTb\endcsname%
      \put(1078,704){\makebox(0,0)[r]{\strut{} 0.2}}%
      \put(1078,1213){\makebox(0,0)[r]{\strut{} 0.3}}%
      \put(1078,1722){\makebox(0,0)[r]{\strut{} 0.4}}%
      \put(1078,2231){\makebox(0,0)[r]{\strut{} 0.5}}%
      \put(1078,2740){\makebox(0,0)[r]{\strut{} 0.6}}%
      \put(1078,3248){\makebox(0,0)[r]{\strut{} 0.7}}%
      \put(1078,3757){\makebox(0,0)[r]{\strut{} 0.8}}%
      \put(1078,4266){\makebox(0,0)[r]{\strut{} 0.9}}%
      \put(1078,4775){\makebox(0,0)[r]{\strut{} 1}}%
      \put(1210,484){\makebox(0,0){\strut{}-0.04}}%
      \put(2153,484){\makebox(0,0){\strut{}-0.03}}%
      \put(3096,484){\makebox(0,0){\strut{}-0.02}}%
      \put(4040,484){\makebox(0,0){\strut{}-0.01}}%
      \put(4983,484){\makebox(0,0){\strut{} 0}}%
      \put(5926,484){\makebox(0,0){\strut{} 0.01}}%
      \put(6869,484){\makebox(0,0){\strut{} 0.02}}%
      \put(308,2739){\rotatebox{-270}{\makebox(0,0){\strut{}\Large $\theta_1 - \dot{\theta}_1$}}}%
      \put(4039,154){\makebox(0,0){\strut{}\Large $\theta_1 + \dot{\theta}_1$}}%
    }%
    \gplgaddtomacro\gplfronttext{%
    }%
    \gplbacktext
    \put(1240,700){\includegraphics{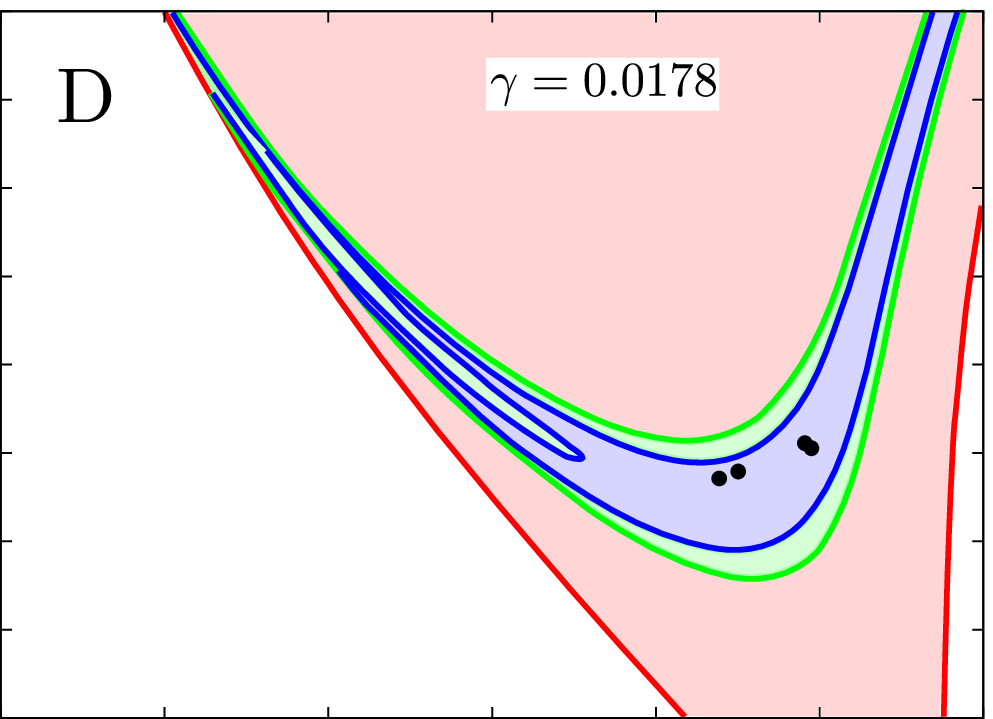}}%
    \gplfronttext
  \end{picture}%
\endgroup

%% file: 0_01903_basin.tex
\begingroup
  \makeatletter
  \providecommand\color[2][]{%
    \GenericError{(gnuplot) \space\space\space\@spaces}{%
      Package color not loaded in conjunction with
      terminal option `colourtext'%
    }{See the gnuplot documentation for explanation.%
    }{Either use 'blacktext' in gnuplot or load the package
      color.sty in LaTeX.}%
    \renewcommand\color[2][]{}%
  }%
  \providecommand\includegraphics[2][]{%
    \GenericError{(gnuplot) \space\space\space\@spaces}{%
      Package graphicx or graphics not loaded%
    }{See the gnuplot documentation for explanation.%
    }{The gnuplot epslatex terminal needs graphicx.sty or graphics.sty.}%
    \renewcommand\includegraphics[2][]{}%
  }%
  \providecommand\rotatebox[2]{#2}%
  \@ifundefined{ifGPcolor}{%
    \newif\ifGPcolor
    \GPcolortrue
  }{}%
  \@ifundefined{ifGPblacktext}{%
    \newif\ifGPblacktext
    \GPblacktexttrue
  }{}%
  \let\gplgaddtomacro\g@addto@macro
  \gdef\gplbacktext{}%
  \gdef\gplfronttext{}%
  \makeatother
  \ifGPblacktext
    \def\colorrgb#1{}%
    \def\colorgray#1{}%
  \else
    \ifGPcolor
      \def\colorrgb#1{\color[rgb]{#1}}%
      \def\colorgray#1{\color[gray]{#1}}%
      \expandafter\def\csname LTw\endcsname{\color{white}}%
      \expandafter\def\csname LTb\endcsname{\color{black}}%
      \expandafter\def\csname LTa\endcsname{\color{black}}%
      \expandafter\def\csname LT0\endcsname{\color[rgb]{1,0,0}}%
      \expandafter\def\csname LT1\endcsname{\color[rgb]{0,1,0}}%
      \expandafter\def\csname LT2\endcsname{\color[rgb]{0,0,1}}%
      \expandafter\def\csname LT3\endcsname{\color[rgb]{1,0,1}}%
      \expandafter\def\csname LT4\endcsname{\color[rgb]{0,1,1}}%
      \expandafter\def\csname LT5\endcsname{\color[rgb]{1,1,0}}%
      \expandafter\def\csname LT6\endcsname{\color[rgb]{0,0,0}}%
      \expandafter\def\csname LT7\endcsname{\color[rgb]{1,0.3,0}}%
      \expandafter\def\csname LT8\endcsname{\color[rgb]{0.5,0.5,0.5}}%
    \else
      \def\colorrgb#1{\color{black}}%
      \def\colorgray#1{\color[gray]{#1}}%
      \expandafter\def\csname LTw\endcsname{\color{white}}%
      \expandafter\def\csname LTb\endcsname{\color{black}}%
      \expandafter\def\csname LTa\endcsname{\color{black}}%
      \expandafter\def\csname LT0\endcsname{\color{black}}%
      \expandafter\def\csname LT1\endcsname{\color{black}}%
      \expandafter\def\csname LT2\endcsname{\color{black}}%
      \expandafter\def\csname LT3\endcsname{\color{black}}%
      \expandafter\def\csname LT4\endcsname{\color{black}}%
      \expandafter\def\csname LT5\endcsname{\color{black}}%
      \expandafter\def\csname LT6\endcsname{\color{black}}%
      \expandafter\def\csname LT7\endcsname{\color{black}}%
      \expandafter\def\csname LT8\endcsname{\color{black}}%
    \fi
  \fi
  \setlength{\unitlength}{0.0500bp}%
  \begin{picture}(7200.00,5040.00)%
    \gplgaddtomacro\gplbacktext{%
      \csname LTb\endcsname%
      \put(1078,704){\makebox(0,0)[r]{\strut{} 0.2}}%
      \put(1078,1213){\makebox(0,0)[r]{\strut{} 0.3}}%
      \put(1078,1722){\makebox(0,0)[r]{\strut{} 0.4}}%
      \put(1078,2231){\makebox(0,0)[r]{\strut{} 0.5}}%
      \put(1078,2740){\makebox(0,0)[r]{\strut{} 0.6}}%
      \put(1078,3248){\makebox(0,0)[r]{\strut{} 0.7}}%
      \put(1078,3757){\makebox(0,0)[r]{\strut{} 0.8}}%
      \put(1078,4266){\makebox(0,0)[r]{\strut{} 0.9}}%
      \put(1078,4775){\makebox(0,0)[r]{\strut{} 1}}%
      \put(1210,484){\makebox(0,0){\strut{}-0.04}}%
      \put(2153,484){\makebox(0,0){\strut{}-0.03}}%
      \put(3096,484){\makebox(0,0){\strut{}-0.02}}%
      \put(4040,484){\makebox(0,0){\strut{}-0.01}}%
      \put(4983,484){\makebox(0,0){\strut{} 0}}%
      \put(5926,484){\makebox(0,0){\strut{} 0.01}}%
      \put(6869,484){\makebox(0,0){\strut{} 0.02}}%
      \put(308,2739){\rotatebox{-270}{\makebox(0,0){\strut{}{\Large $\theta_1 - \dot{\theta}_1$}}}}%
      \put(4039,154){\makebox(0,0){\strut{}{\Large $\theta_1 + \dot{\theta}_1$}}}%
    }%
    \gplgaddtomacro\gplfronttext{%
    }%
    \gplbacktext
    \put(1240,700){\includegraphics{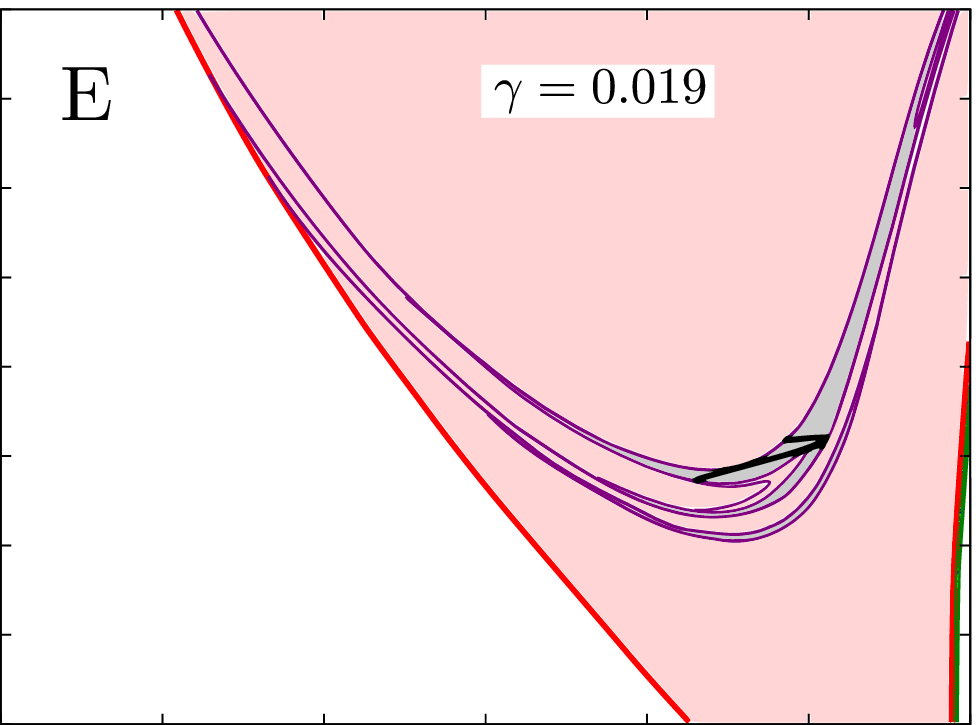}}%
    \gplfronttext
  \end{picture}%
\endgroup

%% file: 0_01903.tex
\begingroup
  \makeatletter
  \providecommand\color[2][]{%
    \GenericError{(gnuplot) \space\space\space\@spaces}{%
      Package color not loaded in conjunction with
      terminal option `colourtext'%
    }{See the gnuplot documentation for explanation.%
    }{Either use 'blacktext' in gnuplot or load the package
      color.sty in LaTeX.}%
    \renewcommand\color[2][]{}%
  }%
  \providecommand\includegraphics[2][]{%
    \GenericError{(gnuplot) \space\space\space\@spaces}{%
      Package graphicx or graphics not loaded%
    }{See the gnuplot documentation for explanation.%
    }{The gnuplot epslatex terminal needs graphicx.sty or graphics.sty.}%
    \renewcommand\includegraphics[2][]{}%
  }%
  \providecommand\rotatebox[2]{#2}%
  \@ifundefined{ifGPcolor}{%
    \newif\ifGPcolor
    \GPcolortrue
  }{}%
  \@ifundefined{ifGPblacktext}{%
    \newif\ifGPblacktext
    \GPblacktexttrue
  }{}%
  \let\gplgaddtomacro\g@addto@macro
  \gdef\gplbacktext{}%
  \gdef\gplfronttext{}%
  \makeatother
  \ifGPblacktext
    \def\colorrgb#1{}%
    \def\colorgray#1{}%
  \else
    \ifGPcolor
      \def\colorrgb#1{\color[rgb]{#1}}%
      \def\colorgray#1{\color[gray]{#1}}%
      \expandafter\def\csname LTw\endcsname{\color{white}}%
      \expandafter\def\csname LTb\endcsname{\color{black}}%
      \expandafter\def\csname LTa\endcsname{\color{black}}%
      \expandafter\def\csname LT0\endcsname{\color[rgb]{1,0,0}}%
      \expandafter\def\csname LT1\endcsname{\color[rgb]{0,1,0}}%
      \expandafter\def\csname LT2\endcsname{\color[rgb]{0,0,1}}%
      \expandafter\def\csname LT3\endcsname{\color[rgb]{1,0,1}}%
      \expandafter\def\csname LT4\endcsname{\color[rgb]{0,1,1}}%
      \expandafter\def\csname LT5\endcsname{\color[rgb]{1,1,0}}%
      \expandafter\def\csname LT6\endcsname{\color[rgb]{0,0,0}}%
      \expandafter\def\csname LT7\endcsname{\color[rgb]{1,0.3,0}}%
      \expandafter\def\csname LT8\endcsname{\color[rgb]{0.5,0.5,0.5}}%
    \else
      \def\colorrgb#1{\color{black}}%
      \def\colorgray#1{\color[gray]{#1}}%
      \expandafter\def\csname LTw\endcsname{\color{white}}%
      \expandafter\def\csname LTb\endcsname{\color{black}}%
      \expandafter\def\csname LTa\endcsname{\color{black}}%
      \expandafter\def\csname LT0\endcsname{\color{black}}%
      \expandafter\def\csname LT1\endcsname{\color{black}}%
      \expandafter\def\csname LT2\endcsname{\color{black}}%
      \expandafter\def\csname LT3\endcsname{\color{black}}%
      \expandafter\def\csname LT4\endcsname{\color{black}}%
      \expandafter\def\csname LT5\endcsname{\color{black}}%
      \expandafter\def\csname LT6\endcsname{\color{black}}%
      \expandafter\def\csname LT7\endcsname{\color{black}}%
      \expandafter\def\csname LT8\endcsname{\color{black}}%
    \fi
  \fi
  \setlength{\unitlength}{0.0500bp}%
  \begin{picture}(7200.00,5040.00)%
    \gplgaddtomacro\gplbacktext{%
      \csname LTb\endcsname%
      \put(1078,704){\makebox(0,0)[r]{\strut{} 0.2}}%
      \put(1078,1213){\makebox(0,0)[r]{\strut{} 0.3}}%
      \put(1078,1722){\makebox(0,0)[r]{\strut{} 0.4}}%
      \put(1078,2231){\makebox(0,0)[r]{\strut{} 0.5}}%
      \put(1078,2740){\makebox(0,0)[r]{\strut{} 0.6}}%
      \put(1078,3248){\makebox(0,0)[r]{\strut{} 0.7}}%
      \put(1078,3757){\makebox(0,0)[r]{\strut{} 0.8}}%
      \put(1078,4266){\makebox(0,0)[r]{\strut{} 0.9}}%
      \put(1078,4775){\makebox(0,0)[r]{\strut{} 1}}%
      \put(1210,484){\makebox(0,0){\strut{}-0.04}}%
      \put(2153,484){\makebox(0,0){\strut{}-0.03}}%
      \put(3096,484){\makebox(0,0){\strut{}-0.02}}%
      \put(4040,484){\makebox(0,0){\strut{}-0.01}}%
      \put(4983,484){\makebox(0,0){\strut{} 0}}%
      \put(5926,484){\makebox(0,0){\strut{} 0.01}}%
      \put(6869,484){\makebox(0,0){\strut{} 0.02}}%
      \put(308,2739){\rotatebox{-270}{\makebox(0,0){\strut{}{\Large $\theta_1 - \dot{\theta}_1$}}}}%
      \put(4039,154){\makebox(0,0){\strut{}{\Large $\theta_1 + \dot{\theta}_1$}}}%
    }%
    \gplgaddtomacro\gplfronttext{%
    }%
    \gplbacktext
    \put(1240,700){\includegraphics{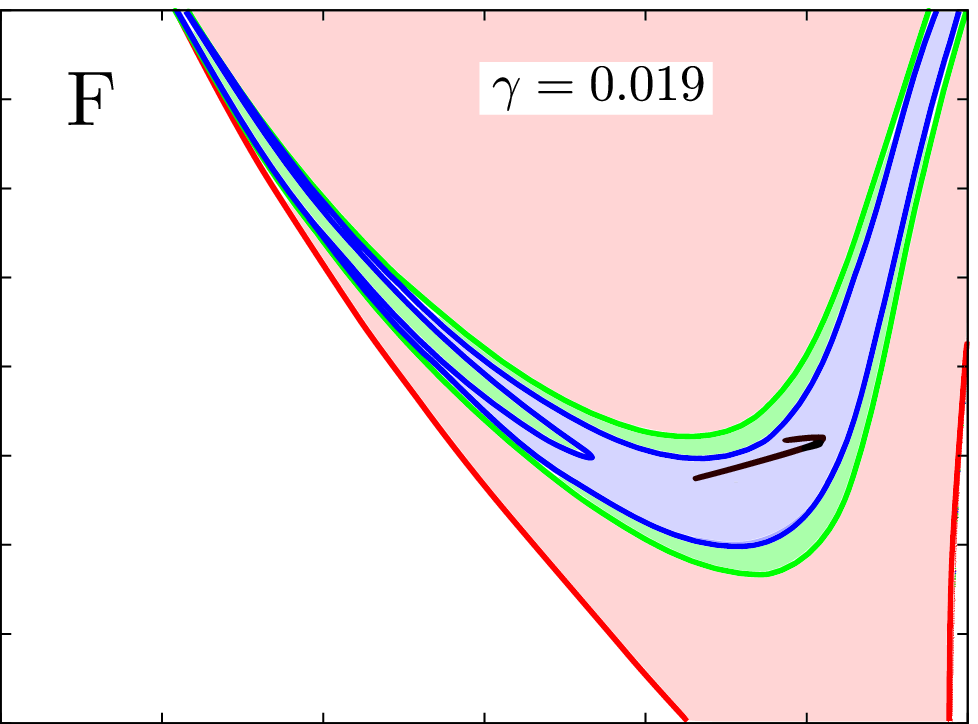}}%
    \gplfronttext
  \end{picture}%
\endgroup

%% file: manuscript.bbl
\begin{thebibliography}{10}

\bibitem{Cavagna1963}
G.~A. Cavagna, F.~P. Saibene, and R.~Margaria, ``External work in walking,''
  {\em Journal of Applied Physiology}, vol.~18, no.~1, pp.~1--9, 1963.

\bibitem{Cavagna1966}
G.~A. Cavagna and R.~Margaria, ``Mechanics of walking.,'' {\em Journal of
  Applied Physiology}, vol.~21, no.~1, pp.~271--278, 1966.

\bibitem{Cavagna1977}
G.~A. Cavagna, N.~C. Heglund, and C.~R. Taylor, ``Mechanical work in
  terrestrial locomotion: two basic mechanisms for minimizing energy
  expenditure,'' {\em American Journal of Physiology - Regulatory, Integrative
  and Comparative Physiology}, vol.~233, no.~5, pp.~R243--R261, 1977.

\bibitem{Ogihara2011}
N.~Ogihara, S.~Aoi, Y.~Sugimoto, K.~Tsuchiya, and M.~Nakatsukasa, ``Forward
  dynamic simulation of bipedal walking in the japanese macaque: Investigation
  of causal relationships among limb kinematics, speed, and energetics of
  bipedal locomotion in a nonhuman primate,'' {\em American Journal of Physical
  Anthropology}, vol.~145, no.~4, pp.~568--580, 2011.

\bibitem{Alexander1980}
R.~Alexander, ``Mechanics of bipedal locomotion,'' in {\em Perspectives in
  Experimental Biology 1} (P.~Spencer-Davies, ed.), pp.~493--504, Oxford:
  Pergamon Press, 1980.

\bibitem{Mochon1980a}
S.~Mochon and T.~A. McMahon, ``Ballistic walking,'' {\em Journal of
  Biomechanics}, vol.~13, no.~1, pp.~49--57, 1980.

\bibitem{Mochon1980b}
S.~Mochon and T.~A. McMahon, ``Ballistic walking: an improved model,'' {\em
  Mathematical Biosciences}, vol.~52, no.~3-4, pp.~241--260, 1980.

\bibitem{Kuo2001}
A.~D. Kuo, ``A simple model of bipedal walking predicts the preferred
  speed-step length relationship,'' {\em ASME Journal of Biomechanical
  Engineering}, vol.~123, no.~3, pp.~264--269, 2001.

\bibitem{Kuo2001b}
A.~D. Kuo, ``Energetics of actively powered locomotion using the simplest
  walking model,'' {\em ASME Journal of Biomechanical Engineering}, vol.~124,
  no.~1, pp.~113--120, 2001.

\bibitem{Srinivasan2007}
M.~Srinivasan and A.~Ruina, ``Idealized walking and running gaits minimize
  work,'' {\em Proceedings of the Royal Society of London A: Mathematical,
  Physical and Engineering Sciences}, vol.~463, no.~2086, pp.~2429--2446, 2007.

\bibitem{Macdonald2014}
J.~H. Macdonald, ``Lateral excitation of bridges by balancing pedestrians,''
  {\em Proceedings of the Royal Society of London A: Mathematical, Physical and
  Engineering Sciences}, vol.~465, no.~2104, pp.~1055--1073, 2009.

\bibitem{Fujiki20150542}
S.~Fujiki, S.~Aoi, T.~Funato, N.~Tomita, K.~Senda, and K.~Tsuchiya,
  ``Adaptation mechanism of interlimb coordination in human split-belt
  treadmill walking through learning of foot contact timing: a robotics
  study,'' {\em Journal of The Royal Society Interface}, vol.~12, no.~110,
  2015.

\bibitem{McGeer1990}
T.~McGeer, ``Passive dynamic walking,'' {\em The International Journal of
  Robotics Research}, vol.~9, no.~2, pp.~62--82, 1990.

\bibitem{McGeer1993}
T.~McGeer, ``Dynamics and control of bipedal locomotion,'' {\em Journal of
  Theoretical Biology}, vol.~163, no.~3, pp.~277--314, 1993.

\bibitem{Asano2005}
F.~Asano, Z.-W. Luo, and M.~Yamakita, ``Biped gait generation and control based
  on a unified property of passive dynamic walking,'' {\em IEEE Transactions on
  Robotics}, vol.~21, no.~4, pp.~754--762, 2005.

\bibitem{Bruijn2011}
S.~M. Bruijn, D.~J.~J. Bregman, O.~G. Meijer, P.~J. Beek, and J.~H. van
  Die\"{e}n, ``The validity of stability measures: A modelling approach,'' {\em
  Journal of Biomechanics}, vol.~44, no.~13, pp.~2401--2408, 2011.

\bibitem{Chyou2011}
T.~Chyou, G.~F. Liddell, and M.~G. Paulin, ``An upper-body can improve the
  stability and efficiency of passive dynamic walking,'' {\em Journal of
  Theoretical Biology}, vol.~285, no.~1, pp.~126--135, 2011.

\bibitem{Russell2005}
S.~Russell, P.~Sheth, and K.~P. Granata, ``Virtual slope control of a forward
  dynamic bipedal walker,'' {\em ASME Journal of Biomechanical Engineering},
  vol.~127, no.~1, pp.~114--122, 2005.

\bibitem{Coleman1998}
M.~J. Coleman and A.~Ruina, ``An uncontrolled walking toy that cannot stand
  still,'' {\em Physical Review Letters}, vol.~80, no.~16, pp.~3658--3661,
  1998.

\bibitem{Collins2005}
S.~Collins, A.~Ruina, R.~Tedrake, and M.~Wisse, ``Efficient bipedal robots
  based on passive-dynamic walkers,'' {\em Science}, vol.~307, no.~5712,
  pp.~1082--1085, 2005.

\bibitem{Collins2001}
S.~H. Collins, M.~Wisse, and A.~Ruina, ``A three-dimensional passive-dynamic
  walking robot with two legs and knees,'' {\em The International Journal of
  Robotics Research}, vol.~20, no.~7, pp.~607--615, 2001.

\bibitem{Goswami1998}
A.~Goswami, B.~Thuilot, and B.~Espiau, ``A study of the passive gait of a
  compass-like biped robot: Symmetry and chaos,'' {\em The International
  Journal of Robotics Research}, vol.~17, no.~12, pp.~1282--1301, 1998.

\bibitem{Johnston2012}
T.~R. Johnston and M.~Hubbard, ``Optimization of the visco-elastic parameters
  describing the heel-region of a prosthesis,'' {\em Journal of Theoretical
  Biology}, vol.~311, pp.~1--7, 2012.

\bibitem{Kuo1999}
A.~D. Kuo, ``Stabilization of lateral motion in passive dynamic walking,'' {\em
  The International Journal of Robotics Research}, vol.~18, no.~9,
  pp.~917--930, 1999.

\bibitem{Kurz2008}
M.~J. Kurz, T.~N. Judkins, C.~Arellano, and M.~Scott-Pandorf, ``A passive
  dynamic walking robot that has a deterministic nonlinear gait,'' {\em Journal
  of Biomechanics}, vol.~41, no.~6, pp.~1310--1316, 2008.

\bibitem{Kwan2007}
M.~Kwan and M.~Hubbard, ``Optimal foot shape for a passive dynamic biped,''
  {\em Journal of Theoretical Biology}, vol.~248, no.~2, pp.~331--339, 2007.

\bibitem{Roos2010}
P.~E. Roos and J.~B. Dingwell, ``Influence of simulated neuromuscular noise on
  movement variability and fall risk in a {3D} dynamic walking model,'' {\em
  Journal of Biomechanics}, vol.~43, no.~15, pp.~2929--2935, 2010.

\bibitem{Su2007}
J.~L.-S. Su and J.~B. Dingwell, ``Dynamic stability of passive dynamic walking
  on an irregular surface,'' {\em ASME Journal of Biomechanical Engineering},
  vol.~129, no.~6, pp.~802--810, 2007.

\bibitem{Garcia1998}
M.~Garcia, A.~Chatterjee, A.~Ruina, and M.~Coleman, ``The simplest walking
  model: Stability, complexity, and scaling,'' {\em ASME Journal of
  Biomechanical Engineering}, vol.~120, no.~2, pp.~281--288, 1998.

\bibitem{Schwab2001}
A.~L. Schwab and M.~Wisse, ``Basin of attraction of the simplest walking
  model,'' in {\em ASME Design Engineering Technical Conferences}, 2001.

\bibitem{deBoer2010}
T.~de~Boer, M.~Wisse, and F.~C.~T. van~der Helm, ``Virtual slope control of a
  forward dynamic bipedal walker,'' {\em ASME Journal of Biomechanical
  Engineering}, vol.~132, no.~7, p.~071012, 2010.

\bibitem{Gritli2012a}
H.~Gritli, N.~Khraief, and S.~Belghith, ``Period-three route to chaos induced
  by a cyclic-fold bifurcation in passive dynamic walking of a compass-gait
  biped robot,'' {\em Communications in Nonlinear Science and Numerical
  Simulation}, vol.~17, no.~11, pp.~4356--4372, 2012.

\bibitem{Gritli2012b}
H.~Gritli, S.~Belghith, and N.~Khraeif, ``Intermittency and interior crisis as
  route to chaos in dynamic walking of two biped robots,'' {\em International
  Journal of Bifurcation and Chaos}, vol.~22, no.~03, p.~1250056, 2012.

\bibitem{LiYang2012}
Q.~Li and X.-S. Yang, ``New walking dynamics in the simplest passive bipedal
  walking model,'' {\em Applied Mathematical Modelling}, vol.~36, no.~11,
  pp.~5262--5271, 2012.

\bibitem{Robinson2008}
C.~Robinson, {\em Dynamical systems: Stability, symbolic dynamics, and chaos}.
\newblock Studies in Advanced Mathematics, Boca Raton, FL: CRC Press, 2008.

\bibitem{Grebogi1983}
C.~Grebogi, E.~Ott, and J.~A. Yorke, ``Crises, sudden changes in chaotic
  attractors, and transient chaos,'' {\em Physica D}, vol.~7, no.~1-3,
  pp.~181--200, 1983.

\bibitem{Suzuki2012}
Y.~Suzuki, T.~Nomura, M.~Casadio, and P.~Morasso, ``Intermittent control with
  ankle, hip, and mixed strategies during quiet standing: A theoretical
  proposal based on a double inverted pendulum model,'' {\em Journal of
  Theoretical Biology}, vol.~310, pp.~55 -- 79, 2012.

\bibitem{Asai2009}
Y.~Asai, Y.~Tasaka, K.~Nomura, T.~Nomura, M.~Casadio, and P.~Morasso, ``A model
  of postural control in quiet standing: Robust compensation of delay-induced
  instability using intermittent activation of feedback control,'' {\em PLoS
  ONE}, vol.~4, p.~e6169, 07 2009.

\bibitem{Funato2016}
T.~Funato, S.~Aoi, N.~Tomita, and K.~Tsuchiya, ``Smooth enlargement of human
  standing sway by instability due to weak reaction floor and noise,'' {\em
  Royal Society Open Science}.
\newblock in press.

\end{thebibliography}
